\documentclass[article,11pt]{amsart}
\usepackage{times}
\usepackage[margin=1.20in]{geometry}
\usepackage{graphicx}
\usepackage{amssymb}
\theoremstyle{plain}
\newtheorem{theorem}{Theorem}[section]

\theoremstyle{definition}

\newtheorem{problem}{Problem}
\newtheorem{conjecture}[problem]{Conjecture}

\theoremstyle{remark}
\newtheorem{remark}[theorem]{Remark}
\newcommand{\suppress}[1]{}

\numberwithin{equation}{section}

\newcommand{\cP}{{\mathcal P}}
\newcommand{\cS}{{\mathcal S}}

\newcommand{\be}{\begin{equation}}
\newcommand{\ee}{\end{equation}}
\newcommand{\bea}{\begin{eqnarray}}
\newcommand{\eea}{\end{eqnarray}}
\newcommand{\bean}{\begin{eqnarray*}}
\newcommand{\eean}{\end{eqnarray*}}

\usepackage{color}

\title [Mathematics, computation, and games] {\Large Three Puzzles on Mathematics, Computation, and Games}

\author[Gil Kalai]{ Gil Kalai \\ Hebrew University of Jerusalem
and Yale University \email{kalai@math.huji.ac.il}}
\thanks{Work
supported in part by ERC advanced grant 320924,
BSF grant 2014290, and NSF grant DMS-1300120.}

\begin{document}

\begin {abstract}
In this lecture I will talk about three mathematical puzzles
involving mathematics and computation
that have preoccupied me over the years.
The first
puzzle is to understand the amazing success of the
simplex algorithm for linear programming.
The second puzzle
is about errors made when votes are counted during
elections. 
The third puzzle is:
are quantum computers possible?

\end {abstract}

\maketitle


\section {Introduction}

The theory of computing and computer science as a whole are  precious resources for mathematicians.
They bring new questions, new profound ideas, and new perspectives on classical mathematical objects,
and  serve as new areas for applications of mathematics and of mathematical reasoning.
In my lecture I will talk about three mathematical puzzles
involving mathematics and computation (and, at times, other fields)
that have preoccupied me over the years.
The connection between mathematics and computing is especially strong
in my field of combinatorics,
and I believe that being able to personally experience the scientific developments described here
over the last three decades may give my description some added value.
For all three puzzles I will try to describe with some detail
both the large picture 
at hand, and zoom in on topics related to my own work.

\subsection *{Puzzle 1: What can explain the success
of the simplex algorithm?}
Linear programming is the problem of
maximizing a linear function $\phi$ subject to
a system of linear inequalities.
The set of solutions for the linear
inequalities is a convex polyhedron $P$. 
The simplex algorithm was developed by George Dantzig. Geometrically
it can be described by moving from one vertex to a
neighboring vertex of $P$ so as to improve the value of the objective function.
The simplex algorithm is one of the most successful mathematical algorithms.
The explanation of this success is an applied, vaguely stated problem, which is
connected with computers.
The problem has strong relations
to the study of convex polytopes, which fascinated mathematicians
from ancient times 
and which served as a starting point for my own research.

If I were required to choose the single most important mathematical explanation
for the success of the simplex algorithm, my choice would point to a
theorem about another algorithm. I would choose
Khachiyan's 1979 theorem
asserting that there is a polynomial-time algorithm for linear programming.
(Or briefly $LP \in {\bf P}$.) Khachiyan's theorem refers to
the ellipsoid method, and the answer is given in the language of
computational complexity, a language  that
did not exist when the question
was originally raised. 

In Section \ref {s:lp} we will discuss  the mathematics of the simplex algorithm,
convex polytopes, and related mathematical objects.
We will concentrate on the study of diameter of graphs of polytopes
and the discovery of randomized subexponential
variants of the simplex algorithm,
I'll mention recent advances: the disproof of the Hirsch conjecture by Santos and the
connection between linear programming and stochastic games leading to
subexponential
lower bounds, discovered by Friedmann, Hansen, and Zwick,
for certain pivot rules.

\subsection *{Puzzle 2: What are methods of election
that are immune to errors in the counting of votes?}
The second puzzle can be seen in the context of
 understanding and planning of electoral methods.
We all remember the sight of vote recount
in Florida in the 2000 US presidential election.
Is the American electoral system, based on electoral votes, inherently
more susceptible to mistakes than the majority system? And what
is the most stable method?
Together with Itai Benjamini and Oded Schramm we investigated
these and similar problems.
We asked the following question: given that there are two
candidates, and each voter chooses at random
and with equal probability (independently) between them, what is
the stability of the outcome, when in the
vote-counting process one percent of the votes is counted incorrectly? (The
mathematical jargon for these errors is
 "noise.") We defined a measure of noise sensitivity of electoral methods
and found that weighted majority
methods are immune to noise, namely, when the probability of
error is small, the chances that the election outcome will be affected
diminish. We also showed that every stable--to--noise method
is "close" (in some mathematical sense)
to a weighted majority method. In later work, O'Donnell,
Oleszkiewicz, and Mossel showed that
the majority system is most stable to noise among all non-dictatorial methods.

\begin {figure}
\centering
\includegraphics[scale=0.3]{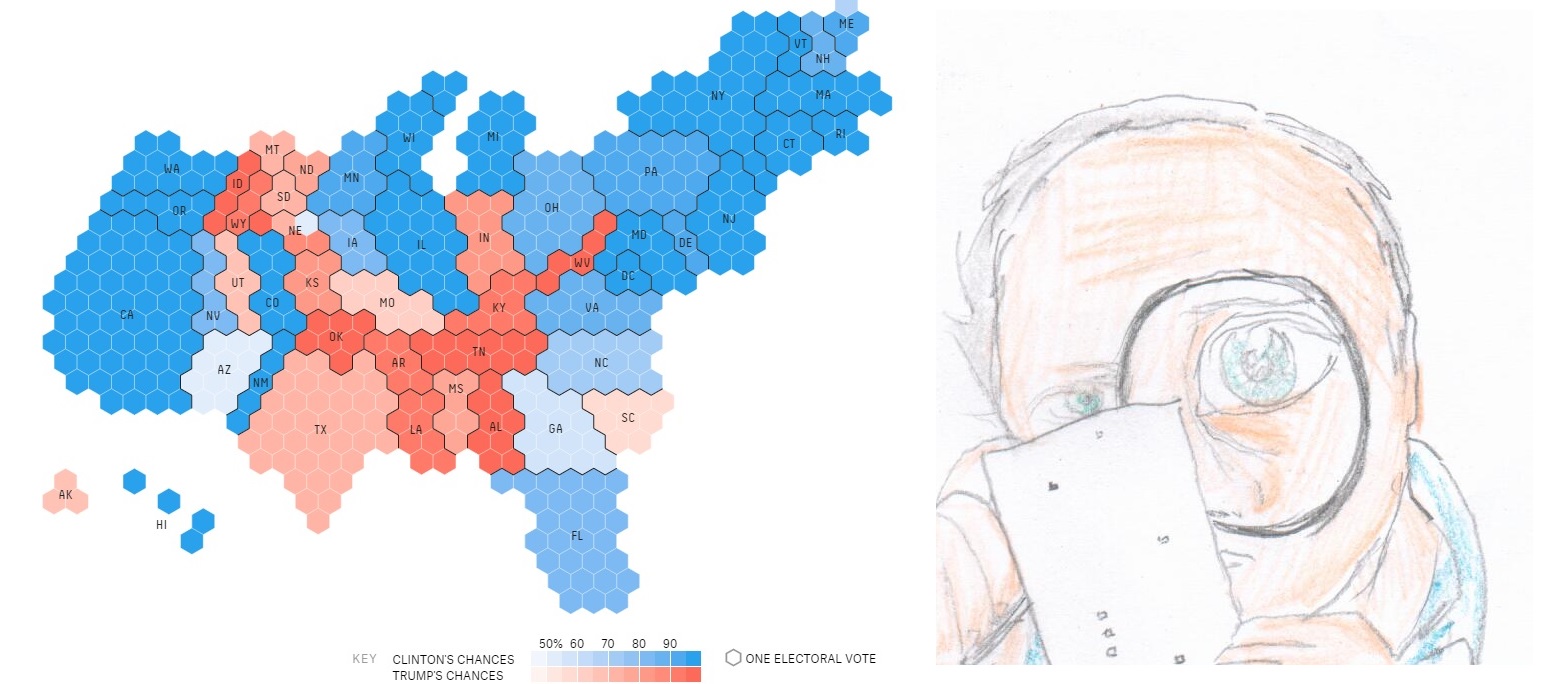}
\caption{\Small{Right: recounts in the 2000 elections (drawing: Neta Kalai).
Left: Hex based demonstration in Nate Silver's site}}
\label{fig:chads}
\end{figure}

Our work was published in 1999, a year before the question appeared in the
headlines in the US presidential election,
and it did not even deal with the subject of elections.
We were interested in understanding the problem of planar percolation, a mathematical model derived from
statistical physics. In our article we
showed that if we adopt an electoral system based
on the model of percolation, this method will be very sensitive to noise. This insight is of no use at
all in planning good electoral methods, but it makes it possible to understand interesting
phenomena in the study of percolation.

After the US presidential election in 2000 we tried to understand the relevance of our
model and the concepts of stability and noise in real-life elections: is the measure for
noise stability  that we proposed relevant, even though the basic assumption that each
voter randomly votes with equal probability for one of the candidates is far from realistic?
The attempt to link mathematical models to questions about
elections (and, more generally, to social science) is fascinating and complicated,
and a pioneer in this study was the Marquis de Condorcet,
a mathematician and philosopher, a democrat, a human rights advocate, and a feminist who
lived in France in the 18th century.
One of Condorcet's findings, often referred to as Condorcet's paradox, is that when there are
three candidates, the majority rule can sometimes lead to cyclic outcomes, and it turns out that the probability
for cyclic outcomes depends on the stability to noise of the voting system.
In Section \ref {s:ns} 
we will discuss noise stability and sensitivity, and various connections
to elections, percolation, games, and computational complexity.

\subsection *{Puzzle 3: Are quantum computers possible?}
A quantum computer is a hypothetical physical device that 
exploits quantum phenomena such as interference
and entanglement in order to enhance computing power.
The study of quantum computation combines fascinating physics,
mathematics, and computer science. In 1994, Peter Shor discovered
that quantum computers would make it possible
to perform certain computational tasks hundreds of orders of magnitude
faster than ordinary computers and, in particular,
would break most of today's encryption methods. At that time, the first doubts about
the model were raised, quantum systems
are of a ``noisy'' and unstable nature. Peter Shor himself found a key to a possible
solution to the problem of ``noise'':
quantum error-correcting codes and quantum fault-tolerance.
In the mid-1990s, three groups of researchers 
showed that ``noisy'' quantum computers still make
it possible to perform all miracles of universal quantum computing,
as long as engineers succeeded in lowering the noise level below a certain threshold.

\begin{figure}
\centering
\includegraphics[scale=0.5]{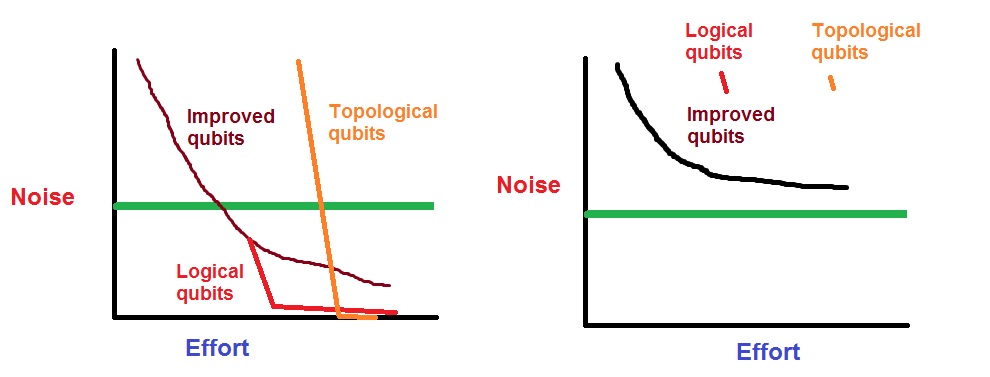}
\caption{{\Small {\bf Left} -- It is commonly believed that by putting more effort for creating qubits
the noise level can be pushed down to as  close to zero as we want.
Once the noise level is small enough and crosses the green 
``threshold'' line,
quantum error correction allows logical qubits to reduce the noise even
further with a small amount of additional effort. Very high quality
topological qubits are also expected.
This belief is supported by ``Schoelkopf's law,'' 
the quantum computing analogue of Moore's law. 
{\bf Right} -- My analysis gives good reasons to expect that we will not be able to reach 
the 
threshold line, that
all attempts for good quality logical and topological qubits will fail, and that 
Schoelkopf's law will break down before 
useful qubits could be created.}}
\label{fig:3}
\end{figure}


A widespread opinion
is that the construction of quantum computers is possible,
that the remaining challenge is essentially of an engineering nature,
and that such computers will be built in the
coming decades. Moreover, people expect to build in the next few years
quantum codes of the quality required for quantum
fault-tolerance,  and to demonstrate the 
concept of "quantum computational supremacy"
on quantum computers with 50 qubits. 
My position
is that it will not be possible to
construct quantum codes that are required for quantum computation,
nor will it be possible to demonstrate quantum computational
superiority in other quantum systems. 
My analysis is based on the same model of noise that led researchers in the 1990s to optimism about
quantum computation, it points to the need for different analyses on different scales,
and it shows that noisy
quantum computers in the small scale (a few dozen qubits) express such a
primitive computational power that it will not
allow the creation of quantum codes that are required as building blocks
for quantum computers on a higher scale.

\subsection *{Near term plans for ``quantum supremacy''}

By the end of 2017\footnote{Of course, for such a major 
scientific project, a delay of a few months and even a couple of years 
is reasonable.}, John Martinis' group is planning to conclude a 
decisive experiment for demonstrating ``quantum supremacy'' on a 50-qubit quantum computer.
(See: Boxio et als (2016) arXiv:1608.00263). 
As they write in the abstract ``A critical question for the field of 
quantum computing in the near future is 
whether quantum devices 
without error correction can perform a well-defined computational task 
beyond the capabilities of state-of-the-art classical computers, 
achieving so-called 
quantum supremacy''
The group intends to study ``the task of sampling from the output distributions 
of (pseudo-)random quantum circuits, 
a natural task for benchmarking quantum computers.'' The objective of this experiment is 
to fix a pseudo-random 
circuit,  run it many times starting from a given initial state to 
create a target state,  and then measure the outcome to reach a probability distribution on 
0-1 sequences of length 50. 

The analysis described in 
Section \ref {s:qc} (based on Kalai and Kindler (2014)) suggests that the outcome 
of this experiment will have 
vanishing correlation with the outcome expected on the ``ideal'' evolution, and that 
the experimental outcomes are actually very very easy to simulate  
classically. They represent distributions that can be expressed by 
low degree polynomials. 
Testing   
quantum supremacy via pseudo-random circuits, against the 
alternative suggested by Kalai, Kindler (2014), 
can be carried out already with a smaller number of qubits (see Fig. 3), 
and
even the 9-qubit experiments (Neil et als (2017)arXiv:1709.06678) should be examined. 
 
The argument for why quantum computers are infeasible is simple.

First, the answer to the question whether quantum devices 
without error correction can perform a well-defined computational task 
beyond the capabilities of state-of-the-art classical computers, is negative. 
The reason is that devices without error correction are computationally very primitive, and 
primitive-based supremacy is not possible. 

Second, the task of creating quantum error-correcting codes 
is harder than the task of demonstrating quantum supremacy,

\begin{figure}
\label{fig:4.01}
\centering
\includegraphics[scale=0.45]{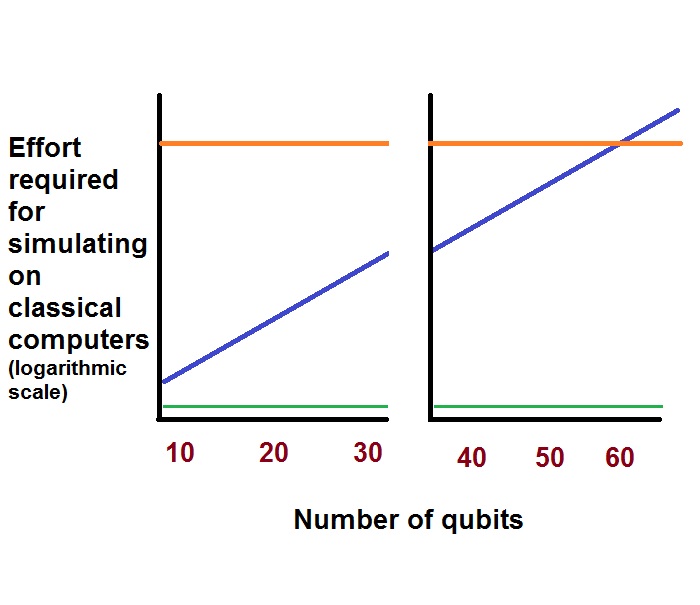}
\caption{\Small{
Quantum supremacy via pseudo-random circuits can be tested for 10-30 qubits. 
The orange line represents the limit for classical computers. 
Kalai, Kindler (2014) suggests close-to-zero correlation 
between the experimental results and outcomes 
based on the noiseless evolution, and further suggests that 
the experimental results are very easy to simulate (the green line).}}
\end{figure}


Quantum computers are discussed in Section \ref {s:qc}, 
we first describe the model, then explain the argument 
for why quantum computers
are not feasible, next we describe predictions for current and near-future devices 
and finally draw some predictions for general noisy quantum systems. 
Section \ref {s:qc} 
presents my 
research since 2005. It is possible, however, that decisive evidence against 
my analysis will be announced or presented in a matter of a few days or a bit later. This is a risk
that I and the reader will have to take.

\subsection * {Perspective and resources}
For books on linear programming see Matou\v{s}ek  and G\"artner (2007)  
and Schrijver (1986). See also Schrijver's (2003) three-volume book on combinatorial optimization,
and a survey article by Todd (2002) on the many facets of linear programming.
For books on convex polytopes see Ziegler's book (1995) 
and Gr\"unbaum's book (1967). 
For game theory, let me recommend the books by
Maschler, Solan and Zamir (2013) 
and Karlin and Peres (2017). 
For books on computational complexity, the reader is referred to Goldreich (2010, 2008),
Arora and Barak (2009), 
and Wigderson (2017, available on the author's homepage) 
For books on Boolean functions and noise sensitivity see O'Donnell (2014) 
and Garban and Steif (2014). 
The discussion in Section \ref {s:ns} complements my 7ECM survey article, 
Boolean functions; Fourier, thresholds and noise. 
It is also related to Kalai and Safra's (2006) survey on threshold phenomena and influence.
For quantum information and computation the reader is referred to
Nielsen and Chuang (2000). 
The discussion in Section \ref{s:qc} follows 
my Notices AMS paper (2016) and its expanded version on the arxive (which is also a good source for references). 
My work has greatly benefited from Internet blog discussions with Aram Harrow, and others, over Regan and Lipton's blog, 
my blog, and other places. 
\begin {remark} These days, references can easily be found through the authors' names and year of publication. 
Therefore, and due also to space limitations, we will provide full references only 
to a few, recently published, papers. 
\end {remark}
\begin {remark} 
Crucial predictions regarding quantum computers
are going to be tested in 
the near future, perhaps even in a few months. I hope to post an updated and more detailed 
version of this paper (with a full bibliography) by the end of 2019.
\end {remark}

\section {Linear programming, polytopes, and the simplex algorithm}
\begin {center}
{\color {blue} 
{\it To Micha A. Perles and Victor L. Klee who educated me as a mathematician.}
}
\end {center}

\label {s:lp}


\subsection {Linear programming and the simplex algorithm}


A linear programming problem is the problem of finding the maximum
of a linear functional (called ``a linear objective function'') $\phi$
on $d$ variables subject to a system of
$n$ inequalities. 

{\bf Maximize} 

$c_1x_1+c_2x_2+\cdots c_dx_d$

{\bf subject to} 

$a_{11}x_1+a_{12}x_2+\cdots + a_{1d}x_d \le b_1$

$a_{21}x_1+a_{22}x_2+\cdots + a_{2d}x_d \le b_2$

...

$a_{n1}x_1+a_{n2}x_2+\cdots + a_{nd}x_d \le b_n$

This can be written briefly as: 
{\bf Maximize} $c^tx$, {\bf subject to} 
$Ax \le b$,
where by convention (only here) vectors are column vectors, 
$x=(x_1,x_2, \dots, x_n)$, $b=(b_1,b_2,\dots,b_n)$, $c=(c_1,c_2,\dots,c_n)$ and $A=(a_{ij})_{1 \le i\le n,1 \le j\le d}$.

The set of solutions to the
inequalities is called the {\it feasible polyhedron} and the
simplex algorithm consists of reaching the optimum
by moving from one vertex to a neighboring vertex.
The precise rule for this move is called ``the pivot rule''.

Here is an example where $n=2d$:

{\bf Maximize} 
$x_1+x_2+\cdots +x_d$, 
{\bf subject to:}
$0 \le x_i \le 1,i=1,2,\dots,d$ 

In this case, the feasible polyhedra is the $d$-dimensional cube. 
Although the number of vertices is exponential, $2^d$, for every pivot 
rule it will take at most $d$ steps to reach the optimal vertex $(1,1,\dots,1)$.

The study 
of linear programming and its major applications in economics was pioneered 
by Kantorovich and Koopmans in the early 1940s.
In the late 1940s George Dantzig realized the importance of linear 
programming for planning problems, and 
introduced the simplex algorithm for solving linear
programming problems. Linear programming and the simplex algorithm are
among the most celebrated applications of mathematics. The question can be 
traced back to a 1827 paper by Fourier. (We will come across Joseph Fourier and 
John von Neumann in every puzzle.)

\subsubsection {Local to global principle}

We describe now two basic properties of linear programming.
\begin {itemize}
\item
If $\phi$ is bounded from above on $P$ then
the maximum of $\phi$ on $P$ is attained at a face of $P$,
in particular there is a vertex $v$ for which the maximum is attained.
If $\phi$ is not bounded from above on $P$ then there is an edge
of $P$ on which $\phi$ is not bounded from above.
\item
A sufficient condition for $v$ to be a vertex of $P$ on which
$\phi$ is maximal is that $v$ is a {\it local maximum}, namely 
$\phi(v) \ge \phi(w)$ for every vertex $w$ which is a
neighbor of $v$.

\end{itemize}

An {\it abstract objective function} (AOF) on a polytope $P$ is an ordering 
of the vertices of $P$ such that every face $F$ of $P$ has 
a unique local maximum.

\subsubsection *{Linear programming duality}

A very important aspect of linear programming is duality.
Linear programming duality associates to an LP problem (given as a 
maximization problem) with $d$ variables and 
$n$ inequalities, a dual LP problem (given as a 
minimization problem) with $n-d$ variables and $n$ inequalities with the same solution.
Given an LP problem, the simplex algorithm for the dual problem 
can be seen as a path-following process on vertices of the 
hyperplane arrangement described by the entire hyperplane 
arrangement described by the $n$ inequalities. 
It moves from one dual-feasible vertex to another, where dual-feasible 
vertex is the optimal vertex to a subset of the inequalities.


\subsection {Overview}

{\it Early empirical experience and expectations.}
The performance of the
simplex algorithm is extremely good in practice. In the early days
of linear programming it was  believed that the common pivot rules
reach the optimum in a number of steps which
is polynomial or perhaps even close to linear
in $d$ and $n$. 
A related conjecture by Hirsch asserted that for 
polyhedra defined by $n$ inequalities in $d$ variables,
there is always a path of length at most $n-d$
between every two vertices. 
We overview some developments regarding linear programming and
the simplex algorithms where by ``explanations'' 
we refer to theoretical results that give some 
theoretical support for the excellent behavior of the simplex algorithm,
while by ``concerns'' we refer to results in the opposite direction.
\begin {enumerate}


\item
{\it The Klee-Minty example and worst case behavior (concern 1).}
Klee and Minty (1972)
found
that one of the most common variants of the simplex
algorithm is exponential in the worst case. In their example, the feasible
polyhedron was combinatorially equivalent to a cube, and all of 
its vertices are actually visited by the algorithm. 
Similar results for other pivot rules were subsequently
found by several authors.

\item
{\it Klee-Walkup counterexample for the Hirsch Conjecture (concern 2).}
Klee and Walkup (1967) 
found an example of an unbounded polyhedron for which the Hirsch Conjecture fails.
They additionally showed that also in the bounded case one cannot realize the Hirsch bound by improving paths.
The Hirsch conjecture for polytopes remained open.
On the positive side, Barnette and Larman gave an 
upper bound for the diameter of graphs of $d$-polytopes with $n$ facets 
which are exponential in $d$ but linear in $n$.

\item
{\it $LP \in$ {\bf P}, via the ellipsoid method (explanation 1).}
In 1979 Khachiyan 
proved that $LP \in P$ namely that there is a polynomial time algorithm
for linear programming. This was a major open problem
ever since the complexity classes P and NP where described in the early 1970s.
Khachiyan's proof was based on
Yudin, Nemirovski and Shor's ellipsoid method, which is not practical for LP.

\item
{\it Amazing consequences.}
Gr\"otschel, Lov\'asz, and Schrijver (1981)
found many theoretical applications of 
the ellipsoid method, well beyond its original scope, and
found polynomial time algorithms for several classes of combinatorial optimization problems. 
In particular they showed that semi-definite programming, the problem of 
maximizing a linear objective function on the set of $m$ by $m$ positive definite 
matrices, is in {\bf P}.

\item
{\it Interior points methods (explanation 2).}
For a few years it seemed like  
there is a tradeoff
between theoretical
worst case behavior and practical behavior. This feeling was shattered
with Karmarkar's 1984 interior point method 
and subsequent
theoretical and practical discoveries.

\item 
{\it Is there a strongly polynomial algorithm for LP? (Concern 3)}
All known polynomial-time algorithms for LP require a number of
arithmetic operations which is polynomial in $d$ and $n$ and linear in $L$,
the number of bits required to represent the input. 
Strongly-polynomial algorithms are algorithms where the
number of arithmetic operations is polynomial in $d$ and $n$ and does not depend on $L$, and
no strongly polynomial algorithm for LP is known.


\item
{\it Average case complexity (explanation 3).}
Borgwart (1982) 
and Smale (1983)
pioneered the study of average case
complexity for linear programming. 
Borgwart showed
polynomial average case behavior for a certain model which
exhibit rotational symmetry. 
In 1983, Haimovich and Adler proved that the average length of 
the shadow boundary path
from the bottom vertex to the top, 
for the regions in an arbitrary arrangement of $n$-hyperplanes 
in $\mathbb R^d$ is at most $d$. 
In 1986, Adler and Megiddo, Adler, Shamir, and Karp, and Todd
proved {\it quadratic} upper
bounds for the simplex algorithm for very general random models
that exhibit certain sign invariance.
All these results are for the shadow boundary rule introduced by Gass and Saaty.



\item
{\it Smoothed complexity (explanation 4).}
Spielman and Teng (2004)
showed  that for the shadow-boundary pivot rule, the average number
of pivot steps required for
a random Gaussian perturbation of variance $\sigma$ of an arbitrary LP
problem is polynomial in $d, n$ and $\sigma^{-1}$. (The dependence on $d$ is at least $d^5$.)
For many, the Spielman-Teng result provides
the best known explanation for the good performance of the simplex algorithm.

\item
{\it LP algorithms in fixed dimension.}
Megiddo (1984) found 
for a fixed value of $d$ a linear time
algorithm for LP problems with $n$ variables. Subsequent simple randomized
algorithms were found by Clarkson (1985,1995),  
Seidel (1991)
and Sharir and Welzl (1992).
Sharir and Welzl 
defined a notion of abstract linear programming problems for which their algorithm applies.

\item
{\it Quasi polynomial bounds for the diameter (explanation 5).}
Kalai (1992) 
and Kalai and Kleitman (1992), 
proved a quasipolynomial upper
bound for the diameter of graphs of $d$-polytopes with $n$ facets.

\item
{\it Sub-exponential pivot rules (explanation 6)}
Kalai (1992) and Matou\v{s}ek,  
Sharir, Welzl (1992) proved
that there are randomized pivot steps which require in expectation
a subexponential number of steps $\exp(K \sqrt {n \log d})$. 
One of those algorithms is the Sharir-Welzl algorithm. 

\item
{\it Subexponential lower bounds for abstract problems (concern 4)}
Matou\v{s}ek 
(1994) 
and Matou\v{s}ek and Szab\'o (2006)
found a subexponential lower
bound for the number of steps required by two basic randomized simplex pivot rules,
for abstract linear programs.

\item
{\it Santos (2012) 
found a counterexample to the Hirsch conjecture (concern 5)}.

\item
{\it The connection with stochastic games.}
Ludwig (1995) 
showed that the subexponential randomized pivot rule can be
applied to the problem posed by Condon of finding the value of certain stochastic games.
For these games this is the best known algorithm.

\item
{\it Subexponential lower bounds for geometric problems (concern 6).}
Building on the connection with stochastic games, 
subexponential lower bounds for genuine LP problems for
several randomized pivot rules  were discovered by
Friedman, Hansen, and Zwick (2011, 2014).

\end {enumerate}

Most of the developments listed  above are in the theoretical
side of the linear programming research and there are also many other 
theoretical aspects. Improving the linear algebra aspects of LP 
algorithms,  and
tailoring the algorithm to specific structural 
and sparsity features of optimization tasks,
are both very important and pose interesting mathematical challenges. 
Also of great
importance are widening the scope of applications, and choosing the
right LP modeling to real-life problems. There is also much 
theoretical and practical work on special families of LP problems.

\subsection {Complexity 1: {\bf P}, {\bf NP}, and $LP$}
The complexity of an algorithmic task is the
number of steps required by a computer program to perform the task.
The complexity is given in terms of the input size, 
and usually refers to the worst case behavior given the input size.
An algorithmic task is
in {\bf P} (called ``polynomial'' or ``efficient'') if there is a
computer program that performs
the task in a number of steps
which is bounded above by a polynomial function of the input size.
(In contrast, an algorithmic task which requires an exponential number of steps in terms of the input size
is referred to as "exponential" or "intractable".)

The notion of a non-deterministic algorithm is one of the most important notions
in the theory of computation.
One way to look at non-deterministic algorithms
is to refer to algorithms where some or all
steps of the algorithm are chosen by an almighty oracle.
Decision problems are algorithmic tasks where the output is either ``yes'' or ``no.''
A decision problem is in {\bf NP} if when the answer is yes, it admits a non-deterministic 
algorithm with a polynomial
number of steps in terms of the input size. In other words,
if for every input for which the answer is ``yes,''
there is an efficient proof demonstrating it, namely, 
a polynomial size proof that a polynomial time algorithm can verify. 
An algorithmic task $A$ is {\bf NP}-hard if a subroutine for solving
$A$ allows solving any problem in {\bf NP}
in a polynomial number of steps. An {\bf NP}-complete
problem is an {\bf NP}-hard problem in {\bf NP}.
The papers by Cook (1971), and Levin (1973) 
introducing {\bf P}, {\bf NP}, and {\bf NP}-complete problems, and  
raising the conjecture that {\bf P} $\ne$ {\bf NP}, 
and the paper by Karp (1972)
identifying 21 central algorithmic problems as {\bf NP}-complete, 
are among the scientific highlights of the 20th century.


Graph algorithms play an important role in computational complexity. 
\textsc{Perfect matching}, 
the agorithmic problem of  
deciding if a given graph $G$ contains a perfect matching 
is in {\bf NP} because exhibiting a perfect matching gives 
an efficient proof that a perfect matching exists.
\textsc{Perfect matching} is in {\bf co-NP} (namely, ``not having a perfect matching'' is in {\bf NP}) 
because by a theorem of Tutte, 
if $G$ does not contain a perfect matching
there is a simple efficient way to demonstrate a proof. An algorithm by Edmonds shows that 
\textsc{Perfect matching} is in {\bf P}.
\textsc {Hamiltonian cycle}, the problem of deciding if $G$ contains a Hamiltonian cycle is also in {\bf NP} -- exhibiting a 
Hamiltonian cycle gives an efficient proof for its existence. However this problem is {\bf NP}-complete.

\begin{figure}
\centering
\includegraphics[scale=0.5]{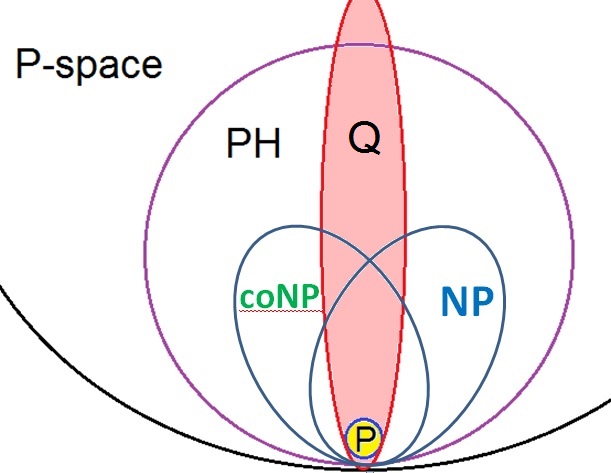}
\caption{\small {The (conjectured) view of some main computational 
complexity classes. The red ellipse represents
efficient quantum algorithms. (See, Section \ref{s:qc}.)}}
\label{fig:0}
\end{figure}

\begin {remark} {\bf P}, {\bf NP}, and {\bf co-NP} are three of the lowest
computational complexity classes in
the {\em polynomial hierarchy {\bf PH}}, which is a countable sequence of such
classes, and there is a rich theory of complexity classes beyond {\bf PH}.
Our understanding of the world of computational complexity depends on a
whole array of conjectures: {\bf NP} $\ne$ {\bf P} is the most famous one. 
A stronger conjecture asserts that {\bf PH} does not collapse,
namely, that there is a strict inclusion between the computational
complexity classes defining the polynomial hierarchy. Counting the 
number of perfect matchings in a graph represents an important complexity class {\bf \#P}
which is beyond the entire polynomial hierarchy.   
\end {remark}

\subsubsection {The complexity of LP, Khachiyan's theorem, and
the quest for strongly polynomial algorithms}

It is known that 
general LP problems can be reduced to the decision problem to
decide if a system of inequalities has a solution. It is therefore easy to see that LP is in {\bf NP}.
All we need is to identify a solution. The duality of linear programming implies that LP is in {\bf co-NP} 
(namely, ``not having a solution'' is in {\bf NP}).
For an LP problem, let $L$ be the number of bits 
required to describe the problem. (Namely, the entries in the 
matrix $A$ and vectors $b$ and $c$.) 

\begin {theorem}[Khachiyan (1979)] 
$LP \in {\bf P}$. The ellipsoid method requires a number of arithmetic steps 
which is polynomial in $n$, $d$, and $L$.
\end {theorem}

The dependence on $L$ in Khachiyan's theorem is linear
and it was met with some surprise. We note that the efficient 
demonstration that a system of linear inequalities has a feasible solution requires 
a number of arithmetic operations which is polynomial in $d$ and $n$ but 
do not depend on $L$. The same applies to an efficient demonstration 
that a system of linear inequalities is infeasible. Also the simplex 
algorithm itself requires a number of arithmetic operations 
that, while not polynomial in $d$ and $n$ in the worst case, does not depend on $L$.
An outstanding open problem is:

\begin {problem}
Is there an algorithm for LP that requires a number of arithmetic
operations which is polynomial in $n$ and $d$ and does not depend on $L$.
\end {problem}

Such an algorithm is called a strongly polynomial algorithm, and this 
problem is one of Smale's ``problems for the 21st century.'' 
Strongly polynomial algorithms are known for various LP problems. 
The Edmonds-Karp algorithm (1972) 
is a strongly polynomial algorithm
for the maximal flow problem. Tardos (1986) 
proved that fixing the feasible
polyhedron (and even fixing only the matrix $A$ used to define the inequalities)  there is a strongly
polynomial algorithm independent of the objective function (and the vector $b$).

\subsection {Diameter of graphs of polytopes and related objects}

\subsubsection {Quasi-polynomial monotone paths to the top}

A few important definitions: a $d$-dimensional polyhedron $P$ is simple if every vertex belongs to $d$ 
edges (equivalently, to $d$ facets.) A linear objective function $\phi$ is generic if 
$\phi(u) \ne \phi (v)$ for two vertices $v \ne  u$. 
The top of $P$ is a vertex for which $\phi$ attains the maximum 
or an edge on which $\phi$ is not bounded. Given a vertex $v$ of $P$ a facet $F$ is active w.r.t. $v$ 
if $\sup_{x \in F} \phi(x) \ge \phi (v)$. 

\begin {theorem}
Let $P$ be a $d$-dimensional simple polyhedron,
let $\phi$ be a generic linear objective function, and let $v$ be a
vertex of $P$. Suppose that there are $n$ active facets w.r.t. $v$.
Then there is a monotone path of length 
$\le n^{\log d+1}$ from $v$ to the top.
\end {theorem}

{Proof:} Let $f(d,n)$ denote the maximum value of the minimum length 
of a monotone path from $v$ to the top. 
(Here, ``the top'' refers to 
either the top vertex or a ray on which $\phi$ is unbounded.)

{\bf Claim:} Starting from a vertex $v$, in $f(d,k)$ steps one can 
reach either the top or vertices in at least $k+1$ active facets.

{\bf Proof:} Let $S$ be a set of $n-k$ active facets. Remove the inequalities defined by these facets to 
obtain a new feasible polyhedron $Q$. If $v$ is not a 
vertex anymore than $v$ belongs to some facet in $S$. 
If $v$ is still a vertex there is a monotone path of length at most $f(d,k)$ from $v$ to the top. 
If one of the edges in the path leaves $P$ then
it reaches a vertex belonging to a facet in $S$. Otherwise it 
reaches the top. Now if a monotone path from $v$ (in $P$) 
of length $f(d,k)$ cannot take us to the top,  there must be such a path that takes us to 
a vertex in some facet belonging to every 
set of $n-k$ facets, and therefore to vertices in at least $k+1$ facets.

{\bf Proof of the Theorem:}
Order the active facets of $P$ according to their top vertex. In $f(d,[n/2])$ steps we can reach 
from $v$ either the top, or a vertex in the top $[n/2]$ active facets. 
In $f(d-1,n-1)$ steps we reach the top $w$ of that facet. 
This will leave us with at most $n/2$ active facets w.r.t. $w$, giving

\begin {equation} 
\label {e:kk}
f(d,n) \le 2f(d,[n/2]) + f(d-1,n-1),
\end {equation}

which implies the bound given by the theorem. (In fact, it gives 
$f(d,n) \le n \cdot {{d+\lceil \log_2 n \rceil} \choose {d}}$.) 
. 

\begin {remark} A  monotone path can be regarded as a 
non-deterministic version of the simplex algorithm where the pivot steps are 
chosen by an oracle.
\end {remark}
\begin {remark}
Let me mention a few of the known upper bounds for the diameter of $d$-polytopes with $n$ 
facets in some special families of polytopes. (Asterisk denotes dual description.) 
Provan, Billera (1980), vertex decomoposable*, (Hirsch bound); 
Adiprasito and Benedetti (2014), flag spheres*, (Hirsch); Nadeff (1989), 0-1 $d$-polytopes, ($d$); 
Balinski (1984), transportation polytopes, (Hirsch); Kalai (1991), dual-neighborly polytopes (polynomial); 
Dyer and Frieze (1994), unimodular, (polynomial); Todd (2014), general polytopes, $(n-d)^{\log d}$ (further 
small improvements followed). 
\end{remark}

\subsubsection {Reductions, abstractions and H\"ahnle's conjecture}

Upper bounds for the diameter are attained at {\it simple} $d$-polytopes, namely $d$-polytopes 
where every vertex belongs to exactly $d$ facets. A more general question deals 
with the dual graphs for triangulations of $(d-1)$-spheres with $n$ vertices. All the known upper bounds 
apply for dual graphs of pure normal $(d-1)$ simplicial complexes. Here ``pure'' means that all 
maximal faces have the same dimension and ``normal'' means that all links of 
dimension one or more are connected. An even more general framework was proposed by 
Eisenbrand, H\"ahnle, Razborov, and Rothvoss (2010). 

\begin {problem}
Consider $t$ pairwise-disjoint nonempty families $F_1,F_2,\dots,F_t$  
of degree $d$ monomials with $n$ variables ($x_1,x_2,\dots,x_n$) 
with the following property:
For every $1\le i <j<k\le t$ if $m_i \in F_i$ and $m_k \in F_k$ then 
there is a monomial $m_j \in F_j$, 
such that the greatest common divisor of $m_i$ and $m_k$ divides $m_j$. 
How large can $t$ be?
\end {problem}

A simple argument shows that the maximum denoted by $g(d,n)$ satisfies relation (\ref{e:kk}).

\begin {conjecture}[H\"ahnle (2010)] 
$g(d,n) \le d(n-1)+1$
\end {conjecture}

One example of equality is taking $F_k$ be all 
monomials $x_{i_1}x_{i_2}\cdots x_{i_d}$ with 
$i_1+i_2+\cdots +i_d= k-d+1$, $k=d,d+1,\dots ,kd$. 
Another example of equality is to let $F_k$ be a single monomial 
of the form $x_i^{d-\ell} x_{i+1}^{\ell}$.
($i=\lfloor k/d \rfloor$ and $\ell = k(\mod~d)$.) 



\subsection {Santos' counterexample}


The $d$-step conjecture is a special case of the Hirsch conjecture known 
to be equivalent to the general case. 
It asserts that a $d$-polytope with $2d$ facets has diameter at most $d$. Santos formulated 
the following strengthening of the $d$-step conjecture:
{\bf Santos' Spindle working Conjecture:} 
Let $P$ be a $d$-polytope with two vertices $u$ and $v$ such that every facet
of $P$ contains exactly one of them. (Such a polytope is called a $d$-spindle.) 
Then the graph-distance between $u$ and $v$ (called simply 
the length of the spindle) is at most 
$d$. Santos proved

\begin{figure}
\centering
\includegraphics[scale=0.5]{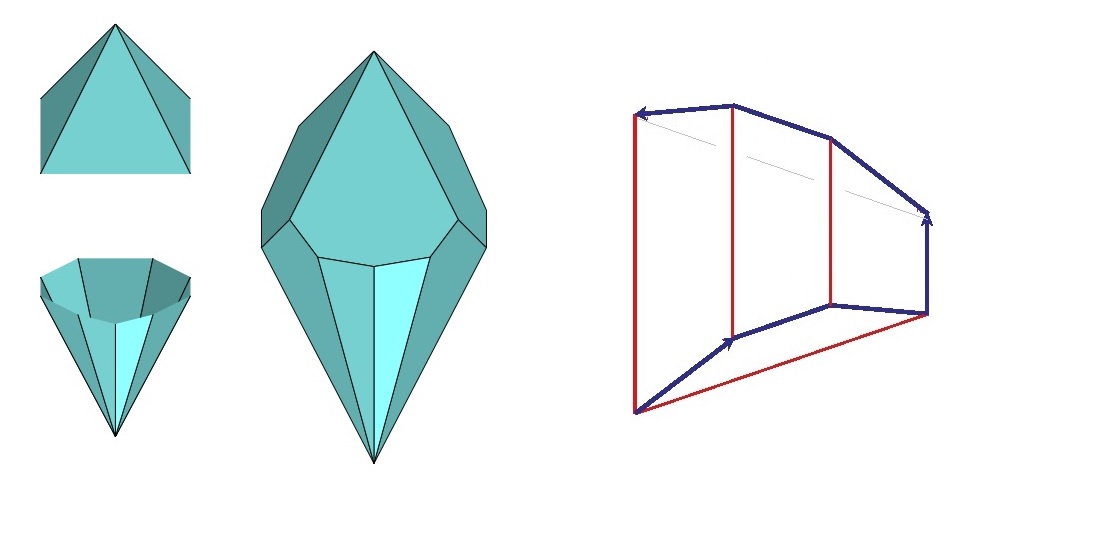}
\caption{{\small Left: a spindle, right: the Klee-Minty cube}}
\label{fig:spin}
\end{figure}


\begin {theorem} [Santos (2014)]

\begin {itemize}
\item [(i)]
The spindle conjecture is equivalent to the Hirsch conjecture. 
More precisely, if there is a $d$-spindle with $n$ facets and length greater than $d$ 
then there is a counter-example to the Hirsch conjecture of 
dimension $n-d$ and with $2n-2d$ facets.
\item[(ii)]
(2) There is a 5-spindle of length 6.
\end {itemize}
\end {theorem}

The initial proof of part (2) had 48 facets and 322 vertices,
leading to a counterexample in dimension 43 with 86 facets and 
estimated to have more than a trillion vertices. 
Matschke, Santos and Weibel (2015) 
found an example with only 25 facets leading 
to a counterexample of the Hirsch conjecture for a 20-polytope with 40 facets and 36,442 vertices. 
An important lesson from Santos' proof is that although reductions are available to 
simple polytopes and simplicial objects, studying the problem for general polytopes has an advantage.
In the linear programming language this tells us that degenerate problems are important. 

\begin {problem}
Find an abstract setting for the diameter problem for 
polytopes which will include graphs of general polytopes, 
dual graphs for normal triangulations, and families of monomials.
\end {problem}




\subsection {Complexity 2: Randomness in algorithms}

One of the most important developments in the theory of computing is the realization 
that adding an internal randomness mechanism can enhance the performance of algorithms.
Some early appearance to this idea came with Monte Carlo methods by 
Ulam, von Neumann and Metropolis, and a factoring algorithm by Berlekamp. 
Since the mid 1970s, and much influenced by Michael Rabin, 
randomized algorithms have become a central paradigm in computer science. 
One of the great achievements were the polynomial time randomized algorithms for testing primality
by Solovay-Strassen (1977) 
and Rabin (1980).  
Rabin's algorithm was related to an earlier breakthrough -- Miller's 
algorithm for primality (1976), 
which was polynomial time conditioned on the validity of the generalized 
Riemann hypothesis. 
The newly randomized algorithms for testing primality were not only theoretically
efficient but also practically excellent! 
Rabin's paper thus gave  ``probabilistic proofs'' that certain large numbers, 
like $2^{400}-593$, are primes,
and this was a new kind of a mathematical proof.
(A deterministic polynomial algorithm for primality 
was achieved  by Agrawal, Kayal, and Saxena (2004).) 
Lov\'asz (1979) 
offered a randomized efficient algorithm for perfect matching in bipartite graphs: Associate to a
bipartite graph $G$ with $n$ vertices on each side, its generic $n \times n$ 
adjacency matrix $A$, where $a_{ij}$ is zero
if the $i$th vertex on one side 
is not adjacent to the $j$th vertex on the other side, 
and $a_{ij}$ is a variable otherwise. Note the
determinant of $A$ is zero if and only if $G$ has no perfect matching. This can be verified
with high probability by replacing $a_{ij}$ with random mod $p$ elements for a large prime $p$.

We have ample empirical experience and some theoretical support to the fact
that pseudo-random number generators are practically sufficient for randomized algorithms.
We also have strong theoretical supports that weak and imperfect sources of randomness are sufficient
for randomized algorithms.


A class of randomized algorithms which are closer to early Monte Carlo algorithms and also to 
randomized algorithms for linear programming, are algorithms
based on random walks. 
Here are two examples:
counting the number of perfect matching for a
general graph $G$ is a ${\bf \#P}$-complete problem. 
Jerrum and Sinclair (1989) and Jerrum, Sinclair and Vigoda (2001) 
found an efficient random-walk based algorithm for
estimating the number of perfect matching up to a multiplicative constant $1+\epsilon$. 
Dyer, Frieze, and Kannan (1991) 
found an efficient
algorithm based on random walks to estimate the volume of a convex body in $\mathbb R^d$. 
Both these algorithms rely in a strong way on the ability to prove a 
spectral gap (or ``expansion'') for various Markov chains. 
Approximate sampling is an important subroutine in the algorithms we have just mentioned and we can regard 
exact or approximate sampling as an important algorithmic task on its own, as
the ability to sample is theoretically and practically important. We mention algorithms by Aldous-Broder (1989) and Wilson (1996) for 
sampling spanning trees 
and an algorithm by 
Randall and Wilson (1999) 
for sampling configurations of the Ising models.

\begin {remark}
The {\it probabilistic method --} 
applying probabilistic methods even for problems with no mention of  probability,
led to major developments in several other mathematical disciplines. In the area of combinatorics
the probabilistic method is especially powerful. (See the book by 
Alon and Spencer (1992).) 
It is an interesting question to what extent proofs obtained by the probabilistic 
method can be transformed into  efficient randomized algorithms.
\end {remark}

\subsection {Subexponential randomized simplex algorithms}

We start with the following simple observation:
Consider the following two sequences. 
The first sequence is defined by $a_1=1$ and $a_{n+1}=a_n+a_{n/2}$, and the second
sequence is defined by $b_1=1$ and $b_{n+1}=b_n+(b_1+\cdots b_n)/n$. Then
$a_n=n^{\theta (\log n)},~~ {\rm and}~~ b_n= e^{\theta (\sqrt n)}.$ 


Next, we describe two basic randomized algorithms for linear programming.

\textsc{Random Edge}: Choose an improving edge uniformly at random and repeat.

\textsc{Random Facet}: 
Given a vertex $v$, choose a facet $F$ containing $v$ uniformly at random. 
Use the algorithm recursively inside $F$ reaching its top $w$, and repeat. 
(When $d=1$, move to the top vertex.)



\textsc{Random Facet} 
(along with some variants) is the first strongly subexponential 
algorithm for linear programming, as well as the first subexponential pivot rule for the simplex algorithm. 


\begin {theorem}[Kalai (1992), Matousek, Sharir, Welzl (1992)]
Let $P$ be a $d$-dimensional simple polyhedron,
let $\phi$ be a linear objective function which is not constant on any edge of $P$, and let $v$ be a
vertex of $P$. Suppose that there are $n$ active facets w.r.t. $v$.
Then \textsc {Random Facet} requires an expected number of at most 
\begin {equation}
\label {e:rf}
\exp (K \cdot \sqrt {n \log d})
\end {equation}
 steps from $v$ to the TOP.
\end {theorem}

{\bf Proof: }
Write $g(d,n)$ for the expected number of pivot steps. The expected 
number of pivot step to reach $w$, the top of the random facet chosen first, is bounded above by $g(d-1,n-1)$.
With probability  $1/d$, $w$ is the $i$th lowest top among 
top vertices of the active facets containing $v$. 
This gives

$$g(d,n) \le g(d-1,n-1) + {\frac {1}{d-1}} \sum_{i=1}^{d-1} g(d,n-i).$$

(Here, we took into account that $v$ itself might be the lowest top.) 
This recurrence relation leads (with some effort) to equation (\ref{e:rf}). $\square$

Note that the argument applies to abstract objective functions on polyhedra. 
(And even, in greater generality, to abstract LP problems as defined by Sharir-Welzl 
and to duals of shellable spheres.)
The appearance of $\exp ({\sqrt n})$ is related to our observation on the sequence $b_n$. 
As a matter of fact, for  the number of steps $G(d)$ for abstract objective functions in the discrete $d$-sphere 
we get the recurrence relation $G(d+1)=G(d)+ (G(1)+G(2)+\cdots +G(d))/d$.
There are few versions of \textsc{Random Facet} that were analyzed (giving slightly lower or 
better upper bound). For the best known 
one see Hansen and Zwick (2015). 
There are a few ideas for improved versions: we can talk about a random face rather than a random facet,
to randomly walk up and down before setting a new threshold, and to try to learn about the problem 
and improve the random choices. The powerful results about lower bounds suggest cautious pessimism. 

\begin {remark} 
Amenta (1994) used Sharir and Welzl's abstract LP problem to settle a Helly type 
conjecture of Gr\"unbaum and Motzkin.
Halman (2004) Considered new large classes of abstract LP problems, found many examples, 
and also related it to Helly-type theorems. 
\end {remark}

\subsubsection {Lower bounds for abstract problems}

As we will see the hope for better upper bounds for \textsc{Random Facet} and related randomized 
versions of the simplex algorithm were faced with formidable examples to the contrary.

\begin {theorem} [Matou\v{s}ek (1994)] 
There exists an abstract objective function  on the $d$-cube on 
which \textsc{Random Facet}  requires on expectation at least $\exp (C \sqrt d)$ steps.
\end {theorem}

Matou\v{s}ek describes a large class of AOF's and showed 
his lower bound to hold in expectation for a randomly chosen AOF. 
G\"artner proved (2002) that for geometric AOF in this 
family, \textsc{Random Facet} requires an expected quadratic time. 

\begin {theorem} [Matou\v{s}ek and Szab\'o (2006)] 
\label {t:matsza}
There exists AOF on the $d$-cube on which \textsc{Random Edge} requires on 
expectation at least $\exp (C d^{1/3})$ steps.
(Hansen and Zwick (2016) improved the bound in Theorem \ref {t:matsza} to $\exp (\sqrt {d \log d})$.) 
\end {theorem}

\subsection {Games 1: Stochastic games, their complexity, and linear programming}



\subsubsection {The complexity of Chess and Backgammon}

Is there a polynomial time algorithm for chess?  Well, if we consider 
the complexity of chess in terms of the board size 
then ``generalized chess'' is very hard. It is {\bf P-space}-complete. 
But if we wish to consider the complexity in terms of 
the number of all possible positions (which for ``generalized chess'' is exponential in the board size), 
given an initial position, it is easy to 
walk on the tree of positions and determine, in a linear number of steps, 
the value of the game.  
(Real life chess is probably intractable, but we note that Checkers was solved.)


Now, what about backgammon?  This question represents one of the most fundamental open 
problems in algorithmic game theory.  
The difference between backgammon and chess is the element of luck; 
in each position your possible moves are determined by a roll of two dice. 

\begin {remark}
Chess and backgammon are games with perfect imformation and their value is achieved by pure strategies.
One of the fundamental insights of game theory is that for zero-sum games with 
imperfect imformation, optimal 
strategies are mixed, namely they are described by random choice between pure strategies. 
For mixed strategies, von Neumann's 1928 minmax theorem asserts that 
a zero sum game with imperfect imformation has a value. 
An optimal strategy for rock-paper-scissors game 
is to play each strategy with equal probability of 1/3. An optimal strategy for 
two-player poker (heads on poker) 
is probably much harder to find. 
\end {remark}

\subsubsection {Stochastic games and Ludwig's theorem}

A simple stochastic game is a two player zero-sum game with perfect 
information, described as follows: We are given 
one shared token and a directed 
graph with two sink vertices labeled '1' and '2' which represent winning 
positions for the two players respectively. All other 
vertices have outdegree 2 and are labelled either by a name of a player 
or as ``neutral''. In addition, one vertex 
is the start vertex. Once the token is on a vertex, the player with the vertex labelling moves, 
and if the vertex is neutral then the move is determined by a toss of a fair coin. 
Backgammon is roughly a game of 
this type. (The vertices represent the player whose turn is to play and the 
outcome of the two dice, and there are neutral vertices representing the rolling of 
the dice. The outdegrees 
are larger than two but this does not make a difference.) This class of 
games was introduced by Condon in 1992.
If there is only one player, 
the game turns into a one-player game with uncertainty, which is  
called a Markov decision process. For Markov 
decision processes finding the optimal strategy is a linear programming problem.

\begin {theorem} [Ludwig (1995)]
There is a subexponential algorithm for solving simple stochastic games
\end {theorem}

The basic idea of the proof is the following: once the strategy of  player 2 is determined 
the game turns into a Markov decision process and the 
optimal strategy for player 1 is determined by solving an LP problem. 
Player one has an optimization problem over the discrete 
cube whose vertices represent his choices in each
vertex labeled by '1'. The crucial observation is that this 
optimization problem defines an abstract objective function 
(with possible equalities) and therefore we can apply \textsc{Random Facet}. 

A more general model of stochastic games with imperfect 
information was introduced by Shapley in 1953. There 
at each step the two players choose actions independently 
from a set of possibilities and their choices determine a reward
and a probability distribution for the next state of the game. 

\begin {problem}[That I learned from Peter Bro Miltersen\\]

\begin {itemize}
\item
{\it Think about backgammon.} Is there a polynomial-time algorithm for finding the value 
of simple stochastic games. 
\item
Can the problem of finding the sink in a unique sink acyclic orientation 
of the $d$-cube be reduced to finding the value of a simple stochastic game?
\item
(Moving away from zero-sum games.) Is there a polynomial-time algorithm (or at least, subexponential algorithm) 
for finding a Nash equilibrium point for a stochastic two-player game (with perfect imformation)?
What about stochastic  
games with perfect imformation with a fixed number of players?
\item
{\it Think about two-player poker}
Is there a polynomial-time algorithm (or at least, subexponential) 
for finding the value of a stochastic zero sum game with imperfect imformation?
\end {itemize}
\end {problem}

\begin {remark}
 What is the complexity for finding objects guaranteed by mathematical theorems?
Papadimitriou (1994) 
developed complexity classes and notions of intractability for  
mathematical methods and tricks! 
(Finding an efficiently describable 
object guaranteed by a mathematical theorem cannot be NP-complete (Megiddo (1988).) 
A motivating conjecture that took many years to prove (in a series of remarkable papers) 
is that Nash equilibria is hard
with respect to PPAD, one of the aforementioned Papadimitriou classes. 
\end {remark}

\begin {problem} [Szab\'o and Welzl (2001)]
How does the problem of finding the sink in a unique sink acyclic orientation of the cube, 
and solving an abstract LP problem, 
fit into Papadimitriou's classes?
\end {problem}

\subsection {Lower bounds for geometric LP problems via stochastic games}

In this section I discuss the remarkable works of Friedmann, Hansen, and Zwick (2011, 2014) 
We talked (having the example of backgammon in mind) 
about  two-player stochastic games with perfect imformation. (Uri Zwick prefers to think of those as 
``games with two and a half players'' with nature being a non strategic player rolling the dice.)
The work of Friedmann, Hansen, and Zwick starts by building 2-player parity games on which 
suitable randomized policy-iteration algorithms perform a subexponential number 
of iterations. Those games are then transformed into 1-player Markov Decision 
Processes (or $1\frac{1}{2}$-player in Uri's view) 
which correspond to concrete linear programs. In their 2014 paper 
they showed a concrete LP problem, where the feasible 
polyhedron is combinatorially a product of simplices, on which \textsc{Random Facet}  
takes an expected number of $\exp({\tilde \Theta (d^{1/3})})$ 
steps, and a variant called \textsc {Random Bland} requires an expected number 
of $\exp({\tilde \Theta (\sqrt d)})$ steps.  
The lower bound even applies to linear programming programs that correspond to shortest paths problems,
(i.e., 1-player games, even by Uri's view) that are very easily solved using
other methods (e.g., Dijkstra (1959)). 
A similar more involved argument gives 
an expected of $\exp (\tilde \Theta (d^{1/4}))$ 
steps for \textsc{Random Edge}! 

\begin {remark} Two recent developments: 
Fearnley and Savani (2015) used the connection between games and LP to show that 
it is {\bf PSPACE}-complete to find the solution that is computed by the simplex method using Dantzig's pivot rule. 
Calude et als. (2017) achieved a quasi-polynomial algorithm for parity games! 
So far, it does not seem that 
the algorithm extends to simple stochastic games 
or has  implications for linear programming.
\end {remark}


\subsection {Discussion}

Is our understanding of the success of the simplex algorithm satisfactory? Are there better practical 
algorithms for semidefinite and convex programming?  Is there a polynomial upper bound for 
the diameter of graphs of $d$-polytopes with $n$ facets?
(Or at least some substantial improvements of known upper bounds). 
Is there a strongly polynomial algorithm for LP? Perhaps even a strongly 
polynomial variant of the simplex algorithm? 
What is the complexity of finding a sink of an acyclic unique-sink orientation of the discrete cube? 
Are there other new interesting efficient, or practically good, algorithms for linear programming? 
What is the complexity of stochastic games? Can a theoretical 
explanation be found for other practically successful algorithms?
(Here, SAT solvers for certain classes of SAT problems, 
and deep learning algorithms come to mind.) Are there good practical algorithms 
for estimating the number of matchings in graphs? for computing volumes for high dimensional polytopes?  
We also face the ongoing challenge of using linear programming and optimization as a source 
for deriving further questions and insights into the 
study of convex polytopes, arrangements of hyperplanes, 
geometry and combinatorics.

\section {Elections and noise}

\begin {center}
{\color {blue}
{\it To Nati Linial and Jeff Kahn who influenced me}
}
\end {center}

\label {s:ns} 

\subsection {Games 2: Questions about voting games and social welfare}

\subsubsection {Cooperative games}

A cooperative game (with side payments) is described by a set of $n$ players $N$, 
and a payoff function $v$ which associates to every  subset $S$ (called {\it coalition}) 
of $N$ a real number $v(S)$. We will assume that $v(\emptyset)=0$. 
Cooperative games were introduced by von Neumann and Morgenstern.
A game is monotone if $v(T) \ge v(S)$ when $S \subset T$. 
A voting game is a monotone cooperative game in which $v(S) \in \{0,1\}$. 
If $v(S)=1$ we call $S$ a winning coalition and if $v(S)=0$ then $S$ is a losing coalition.
Voting games represent voting rules for two-candidate elections, the candidates being Anna and Bianca. 
Anna wins if the set of voters that voted for her is a winning coalition. Important voting rules are
the majority rule, where $n$ is odd and the winning coalitions are those with more than $n/2$ 
voters, 
and the dictatorship rule, where the winning coalitions are those 
containing a fixed voter called ``the dictator.'' Voting games are also referred 
to as monotone Boolean functions. 


\subsubsection {How to measure power?}

There are two related measures of power for voting games 
and both are defined in terms of general cooperative games. 
The Banzhaf measure of power for player $i$, $b_i(v)$ (also called the influence of $i$) is 
the expected value of $v(S \cup \{i\}) - v(S)$ taken over all 
coalitions $S$ that do not contain $i$. The Shapley value of player $i$ is defined as follows:
For a random ordering of the players consider the coalition $S$ of players who 
come before $i$ in the ordering. The Shapley value, $s_i(v)$. is the 
expectation over all $n!$ orderings of $v(S \cup \{i\})- v(S)$.
(For voting games, the Shapley value is also called the Shapley-Shubik power index.) 
For voting games if $v(S)=0$ and $v(S \cup \{i\})=1$, we call voter $i$ {\it pivotal} 
with respect to $S$.  

\subsubsection {Aggregation of information} 

For a voting game $v$ and $p, 0\le p \le 1$ 
denote by $\mu_p(v)$ the probability that a random set $S$ of players is a
winning coalition when for every player $v$ the probability
that $v \in S$ is $p$, independently for all players. Condorcet's
Jury theorem asserts that when $p>1/2$, for the sequence $v_n$ of majority games
on $n$ players $\lim_ {n \to \infty} \mu_p(v_n) = 1.$
This result, a direct consequence of the law of large numbers, is
referred to as asymptotically complete aggregation of information.

A voting game is {\it strong} (also called 
{\it neutral}) if a coalition is winning iff its complement is losing.
A voting game is strongly balanced if 
precisely half of the coalitions are winning and it is balanced  if 
$0.1 \le \mu_{1/2}(v)\le 0.9$. A voting game is weakly symmetric if it is invariant 
under a transitive group of permutations of the voters.


\begin {theorem} [Friedgut and Kalai (1996), Kalai (2004)]

(i) Weakly-symmetric balanced voting games aggregate information.

(ii) Balanced voting games aggregate information iff their maximum Shapley value tend to zero.

\end {theorem}

\subsubsection {Friedgut's Junta theorem}

The total influence, $I(v)$, of a balanced voting game is the sum of 
Banzhaf power indices for all players. 
(Note that the sum of Shapley values of all players is one.)
For the majority rule the total influence is the maximum over all voting games and 
$I=\theta (\sqrt n)$. The total influence for dictatorship is one which is 
the minimum for strongly balanced games. A voting game is a $C$-Junta if there is a  
a set $J$, $|J| \le C$ such that $v(S)$ depends only on $S \cap J$.    

\begin {theorem}[Friedgut's Junta theorem (1998)] 
\label {t:junta}

For every $b, \epsilon >0$ there is $C=C(b,\epsilon)$ with the following property:
For every $\epsilon,b >0$  a voting game $v$ with total influence at most $b$ 
is $\epsilon$-close to a $C$-junta $g$. 
(Here, $\epsilon$-close means that for all but a fraction $\epsilon$ 
of sets $S$,  $v(S)=g(S)$.

\end {theorem}

\subsubsection {Sensitivity to noise}

Let $w_1,w_2, \dots, w_n$ be  nonnegative real weights and $T$ be a real number. A 
weighted majority 
is defined by $v(S)=1$ iff $\sum_{i \in S} w_i \ge T$. 

Consider a two-candidate election based on a voting game $v$ where each voter votes for 
one of the two candidates at random, 
with probability 1/2, and these probabilities are independent. 
Let $S$ be the set of voters voting for Anna, who wins the election if $v(S)=1$.
Next consider a scenario where in the vote counting process there is a small probability $t$ 
for a mistake where the vote is miscounted, and assume that these mistakes are 
also statistically independent. The set of voters believed to vote for Anna after the counting is $T$. 
Define $N_t(v)$ as the probability that $v(T) \ne v(S)$.
A family of voting games is called uniformly noise stable 
if for every $\epsilon>0$ there exists $t>0$ such that 
$N_t(v) < \epsilon$. 
A sequence $v_n$ of strong voting games is noise sensitive 
if for every $t>0$ $\lim_{n \to \infty}  N_t(v_n)=1/2.$

\begin {theorem} [Benjamini, Kalai, Schramm (1999)]
\label {t:bks}
For a sequence of balanced voting games $v_n$ each of the following two conditions implies that
$v_n$ is noise sensitive:

(i) The  maximum correlation between $v_n$ and a balanced weighted majority game tends to 0.

(ii) $\lim_{n\to \infty} \sum _i b_i^2(v_n)=0.$
\end {theorem}

\subsubsection {Majority is stablest}
Let $v_n$ be the majority voting games with $n$ players.
In 1989 Sheppard proved that
$\lim_{n \to \infty} N_t(v_n)= \frac {\arccos (1-2t)}{\pi}.$ 


\begin {theorem}[Mossel, O'Donnell, and Oleszkiewicz (2010)]
\label {t:moo}
Let $v_n$ be a sequence of games with diminishing maximal Banzhaf power index. Then
$$N_t(v_n) \ge \frac {\arccos (1-2t)}{\pi} - o(1).$$ 
\end {theorem}

\subsubsection {The influence of malicious counting errors}
Let $S$ be a set of voters. $I_S(v)$ is the probability over sets of voters $T$ 
which are disjoint from $S$ that
$v (S \cup T)=1$ and $v(T)=0.$  

\begin {theorem} [Kahn, Kalai, Linial (1988)] 
\label {t:kkl}
For every balanced 
voting game $v$.

(i) There exists a voter
$k$ such that $$b_k(f) \ge C  \log n/n.$$

(ii) There exists a set $S$ of $a(n) \cdot n/\log n$ voters,
where $a(n)$ tends to infinity with $n$ as slowly as we wish, such that $I_S(v)=1-o(1)$.
\end {theorem}

This result was conjectured by Ben-Or and Linial (1985) who gave
a  ``tribe'' example showing that both parts of the theorem are sharp.
Ajtai and Linial (1993) found a voting game where no set of $o(n/\log^2(n))$ can influence
the outcome of the elections in favor of even one of the candidates.

\subsubsection {``It ain't over 'till it's over'' theorem}

Consider the majority voting game when the number of voters tends to infinity and every voter 
votes for each candidate with equal probability, independently. There exists (tiny) $\delta >0$ with 
the following property: When you count $99$\%  of votes chosen at random, still
with probability tending to one, condition on the votes counted, 
each candidate has probability larger than $\delta$ of winning. 
We refer to this property of the majority function as the (IAOUIO)-property. 
Clearly, dictatorship and Juntas do not have the (IAOUIO)-property.

\begin {theorem} [Mossel, O'Donnell, and Oleszkiewicz (2010)]
Every sequence of voting games with diminishing maximal Banzhaf power index has the (IAOUIO)-property. 
\end{theorem}
 
\subsubsection {Condorcet's paradox and Arrow's theorem}

A generalized social welfare function is a map from $n$ 
voters' order relations on $m$ alternatives, to a complete antisymmetric relation
for the society, satisfying the following two properties. 

(1) If every voter prefers $a$ to $b$ then so is the society. 
(We do not need to assume that this dependence is monotone.)

(2) Society's preference between $a$ and $b$ depends only 
on the individual preferences between these candidates.

A social welfare function is a generalized welfare function such that for every $n$-tuples of order 
relations of the voters, the society 
preferences are acyclic (``rational''). 

\begin {theorem}[Arrow (1950)]
For three or more alternatives, the only social welfare functions are dictatorial. 
\end {theorem}

\begin {theorem} [Kalai (2002), Mossel (2012), Keller (2012)]
For three or more alternatives the only nearly rational generalized social welfare 
functions are nearly dictatorial.
\end {theorem}

\begin {theorem} [Mossel, O'Donnell, and Oleszkiewicz (2010)]
 The majority gives asymptotically ``most rational'' social preferences 
 among generalized social welfare functions based on 
strong voting games with vanishing maximal Banzhaf power. 
\end {theorem}

A choice function is a rule which, based on individual 
ranking of the candidates, gives the winner of the election. 
Manipulation (also called ``non-naive voting'' 
and ``strategic voting'') is a situation 
where given the preferences of other voters, a voter may 
gain by not being truthful about his preferences. 

\begin {theorem} [Gibbard (1977), Satterthwaite (1975)]
\label {t:gs}
Every non dictatorial choice function is manipulable.
\end {theorem}

\begin {theorem} [Friedgut, Kalai, Keller, Nisan (2011), Isaksson, Kindler, Mossel (2012), Mossel, R\'acz (2015)]
Every nearly non-manipulable choice function is nearly dictatorial. 
\end {theorem}

\subsubsection {Indeterminacy and chaos}

Condorcet's paradox asserts that the majority rule may lead to cyclic outcomes for three candidates. 
A stronger result was proved by McGarvey (1953): every asymmetric preference relation on $m$ alternatives is the outcome 
of majority votes between pairs of alternatives for some individual rational 
preferences (namely, acyclic preferences) for a large number of voters.
This property is referred to as {\it indeterminacy}. A stronger property is that when the individual order relations 
are chosen at random, the probability for every asymetric relation is bounded away from zero. This is 
called {\it stochastic indeterminacy}. Finally, complete chaos refers to a situation where   
in the limit all the probabilities for asymmetric preference relations are the 
same -- $2^{-{{m} \choose {2}}}$.

\begin {theorem} [Kalai (2004, 2007)] 

(i)  Generalized social welfare functions based on voting games that aggregate information 
lead to complete indeterminacy. In particular this applies when the maximum Shapley value tends to zero.

(ii)  Generalized social welfare functions based on voting games 
where the maximum Banzhaf value tends to zero 
leads to stochastic indeterminacy.

(iii) Generalized social welfare functions based on noise-sensitive voting games lead to complete chaos. 

\end {theorem}


\subsubsection {Discussion}
{\it Original contexts for some of the results.}
Voting games are also called monotone Boolean functions and some of the results we discussed 
were proved in this context. Aggregation of information is also 
referred to as the sharp threshold phenomenon, which is important 
in the study of random graphs, percolation theory and other areas. 
Theorem \ref {t:kkl} was studied in the context of distributed computing and the question of 
collective coin flipping: procedures allowing $n$ agents to reach a random bit. 
Theorem \ref {t:bks} was studied in the context of critical planar percolation.
Theorem \ref {t:junta} was studied in the context of the combinatorics and probability 
of Boolean functions. 
The majority is stablest theorem was studied both in the context of hardness of 
approximation for the \textsc{Max Cut} problem (see Section \ref{s:cut}), and  
in the context of social choice. 
Arrow's theorem and Theorem \ref {t:gs} 
had immense impact on theoretical economics and political science.
There is a large body of literature with extensions and interpretations of Arrow's theorem, and related 
phenomena were considered by many. Let me mention the more recent study of judgement aggregation,  
and also Peleg's books (1984, 2010) and Balinski's books (1982, 2010) on voting 
methods that attempt to respond to 
the challenge posed by Arrow's theorem.
Most proofs of the results discussed here go through Fourier analysis of 
Boolean functions that we discuss in Section \ref {s:fourier}. 

{\it Little more on cooperative games.}
I did not tell you yet about the most important 
solution concept in cooperative game theory (irrelevant 
to voting games) -- the core. The core of the game is an assignment 
of $v(N)$ to the $n$ players so that the members of 
every coalition $S$ get together at least $v(S)$. Bondareva and Shapley found necessary and sufficient 
conditions for the core to be non empty, closely related to 
linear programming duality. I also did not talk about games without 
side payments. There, $v(S)$ are sets of vectors which 
describe the possible payoffs for the player in $S$ if they go together.
A famous game with no side payment is Nash's bargaining problem 
for two players. Now, you are just one step away from 
one of the deepest and most beautiful results in game theory, 
Scarf's conditions (1967) for non-emptiness of the core.

{\it But what about real-life elections?}
The relevance and interpretation of mathematical modeling and results regarding voting rules, games, 
economics and social science, 
is a fairly involved matter. It is interesting to examine some notions discussed here in the 
light of election polls which are often based on a 
more detailed model. 
Nate Silver's detailed forecasts provide a special opportunity. 
Silver computes the probability of victory for every candidate based on 
running many noisy  simulations which are in turn based on the outcomes of individual polls. 
The data in Silver's forecast contain an estimation for the event ``recount'' which 
essentially measures noise sensitivity,  
and it would be interesting to compare noise sensitivity in this more realistic scenario to the simplistic model of 
i.i.d. voter's behavior. 
Silver also computes certain power indices based on the probability for pivotality, again,  under his model. 

{\it But what about real-life elections (2)?} Robert Aumann remembers a Hebrew University 
math department meeting convened to choose two new members from among four very serious 
candidates.  The chairman, a world-class mathematician, asked Aumann for a voting procedure.  
Aumann referred him to Bezalel Peleg, an expert on social choice and voting methods.  
The method Peleg suggested was adopted, and two candidates were chosen accordingly.  
The next day, the chairman met Aumann and complained that a majority of the department 
opposes the chosen pair, indeed prefers a specific different pair!  Aumann 
replied, yes, but there is another pair of candidates that the majority prefers to yours, 
and yet another pair that the majority prefers to THAT one; and the pair elected 
is preferred by the majority to that last one! Moreover, there is a theorem 
that says that such situations cannot be avoided under any voting rule. 
The chairman was not happy and said dismissively: ``Ohh, you guys and your theorems.''

\subsection {Boolean functions and their Fourier analysis}


We start with the discrete cube $\Omega_n=\{-1,1\}^n$.
A Boolean function is a map $f:\Omega_n\to \{-1,1\}$.

\begin {remark}A Boolean function represents a family of subsets
of $[n]=\{1,2,\dots,n\}$ (also called hypergraph) which are central objects in extremal combinatorics.   
Of course, voting games are monotone Boolean functions. 
We also note that in agreement with Murphy's law, roughly half of the times it is convenient to consider 
additive notation, namely to
regard $\{0,1\}^n$ as the discrete cube and Boolean functions as functions to $\{0,1\}$.
(The translation is $0 \to 1$ and $1 \to -1$.) 
\end {remark}


\subsubsection {Fourier}
\label {s:fourier} 

Every real function $f:\Omega_n\to \mathbb R$ can
be expressed in terms of the Fourier--Walsh basis. We write here and
for the rest of the paper $[n]=\{1,2,\dots,n\}$.

\begin {equation}
\label{e:fourier}
f=\sum \{\hat f(S)W_S:~ S \subset [n]\},
\end {equation}

where the {\it Fourier-Walsh function} $W_S$ is 
simply the monomial $W_S(x_1,x_2,\dots,x_n)= \prod_{i \in S}x_i$.

Note that we have here $2^n$ functions, one for each subset $S$ of $[n]$. 
Let $\mu$ be the uniform probability measure on $\Omega_n$. The functions $W_S$ 
form an orthonormal basis of $\mathbb R ^{\Omega_n}$
with respect to the inner product $$\langle f,g\rangle =\sum_{x \in \Omega_n}\mu(x)f(x)g(x).$$
The coefficients $\hat f(S)=\langle f,W_s\rangle$, $S \subset [n]$,  in (\ref{e:fourier}) are real numbers,
called the {\it Fourier coefficients} of $f$.
Given a real function $f$ on the discrete
cube with Fourier expansion $f=\sum \{\hat f(S)W_S:~ S \subset [n]\},$
the noisy version of $f$, denoted by $T_\rho (f)$ is defined by
$T_\rho(f)=\sum \{ \hat f(S)(\rho)^{|S|} W_S:~ S \subset [n]\}.$

\subsubsection 
{Boolean formulas, Boolean circuits, and projections.}
(Here it is convenient to think about the additive convention.)
Formulas and circuits allow to build complicated Boolean functions
from simple ones and they have crucial importance in
computational complexity. Starting with $n$ variables $x_1,x_2, \dots, x_n$, a {\it literal} is a
variable $x_i$ or its negation $\neg x_i$.
Every Boolean function can be written as a formula in conjunctive normal form,
namely as AND of ORs of literals.  A {\it circuit}
of depth $d$ is defined inductively as follows. A circuit of depth zero is a literal.
A circuit of depth one consists of an OR or AND {\it gate} applied to a set of literals,
a circuit of depth $k$ consists of an OR or AND gate applied to the outputs of circuits of depth $k-1$.
(We can assume that gates in the odd levels are all OR gates and that
the gates of the even levels are all AND gates.)
The size of a circuit is the number of gates.
Formulas are circuits where we allow to use the output
of a gate as the input of only one other gate.
Given a Boolean function $f(x_1,x_2,\dots,x_n,y_1,y_2,\dots,y_m)$ 
we can look at its 
projection (also called trace) $g(x_1,x_2,\dots,x_n)$ 
on the first $n$ variables. $g(x_1,x_2,\dots,x_n)=1$ if there are 
values $a_1,a_2,\dots,a_m)$ (depending on the $x_i$s) such that 
$f(x_1,x_2,\dots,x_n,a_1,a_2,\dots,a_m)=1$. 
Monotone formulas and circuits are those where all literals are variables. (No negation.)

{\it Graph properties.}
A large important family of examples is obtained as follows.
Consider a property $P$ of graphs on $m$ vertices. Let $n=m(m-1)/2$,
associate Boolean variables with the $n$ edges of the complete graph $K_m$, and represent
every subgraph of $K_m$ by a vector in $\Omega_n$.
The property $P$ is now represented by a Boolean function on $\Omega_n$.
We can also start with an arbitrary graph $H$ with $n$ edges and
for every property $P$ of subgraphs of $H$
obtain a Boolean function of $n$ variables based on $P$.

\subsection {Noise sensitivity everywhere (but mainly percolation)}

One thing we learned through the years is that noise sensitivity 
is (probably) a fairly common phenomenon. 
This is already indicated by Theorem \ref {t:bks}. Proving 
noise sensitivity can be difficult. 
I will  talk in this section about results on the critical planar 
percolation model, and conclude with a problem by Benjamini and Brieussel.
I will not be able to review here many other noise-sensitivity results that 
justify the name of the section.

\subsubsection {Critical planar percolation}
The crossing event for planar percolation refers to an $n$ by $n$  
square grid and to the event, when every edge is chosen with probability 1/2, 
that there is a path crossing from the left side to the right side of the square.

\begin {theorem} [Benjamini, Kalai, and Schramm] (1999) 
The crossing event for percolation is sensitive to $1/o(\log n)$ noise.
\end {theorem}

\begin {theorem} [Schramm and Steif (2011)]  

The crossing event for percolation is sensitive to $(n^{-c+o(1)})$ noise, for some $c>0$.
\end {theorem}


\begin {theorem} [Garban, Pete and Schramm (2010, 2013); Amazing!] 
The crossing event for (hex) percolation is sensitive to $(n^{-(3/4)+o(1)})$ noise. 
The spectral distribution 
has a scaling limit and it is supported by 
Cantor-like sets of Hausdorff dimension 3/4.
\end {theorem}

\begin {remark}[Connection to algorithms]
The proof of Schramm and Steif is closely 
related to the model of computation of random decision trees. Decision tree complexity refers 
to a situation where given a  Boolean function we 
would like to find its value by asking as few as possible questions about specific instances.
Random decision trees allow to add randomization in the choice of the next question. 
These relations are explored in O'Donnell, Saks, Schramm, and Servedio (2005) 
and have been very useful in recent works in percolation theory.
\end {remark}

\begin {remark}[Connection to games]
Critical planar percolation is closely related to 
the famous game of Hex. Peres, Schramm, Sheffield, and Wilson (2007)
studied random turn-Hex where a coin-flip determines the identity 
of the next player to play. They found a simple but surprising 
observation that the value of the game when both players play the random-turn game optimally 
is the same as when both players play randomly. (This applies in 
much greater generality.) Richman considered 
such games which are auction-based turn. Namely, the players 
bid on  who will play the next round. 
A surprising, very general analysis (Lazarus, Loeb, Propp, Ullman (1996)) shows that the 
value of the random-turn game is closely related to 
that of the auction-based game! 
Nash famously showed that for ordinary 
Hex, the first player wins but his proof 
gives no clue as to the winning strategy. 
\end {remark}

\subsubsection {Spectral distribution and Pivotal distribution}

Let $f$ be a monotone Boolean function with $n$ variables. We can associate to $f$ two important probability distributions on 
subsets of $\{1,2,\dots, n\}$. 
The spectral distribution of $f$,  
$\cS (f)$ 
gives a set $S$ a probability $\hat f^2(S)$.  
Given $x \in \Omega_n$ the $i$th variable is {\it pivotal} if when we flip the value 
of $x_i$ the value of $f$ is flipped as well. 
The pivotality distribution 
$\cP (f)$ 
gives a set $S$ the probability that 
$S$ is the set of pivotal variables. 
It is known that the first two moments of $\cS$ and $\cP$ agree.

\begin {problem}
Find further connections between 
$\cS (f)$  and $\cP (f)$ for all Boolean functions and for specific classes of Boolean functions. 
\end {problem}

\begin {conjecture}
Let $f$ represents the crossing event in planar percolation. 
Show that 
$H(\cS )(f))= O(I(f))$ and $H(\cP (f))= O(I(f))$. 
(Here $H$ is the entropy function.) 
\end {conjecture}

The first inequality is a special case of the entropy-influence conjecture of Friedgut and Kalai (1996) 
that applies to general Boolean functions. The second inequality is not so general (it fails for the majority). 
The majority function $f$ has an unusual property (that I call ``the anomaly of majority'') that 
$I(f) = \mathbb E(\cS)= c \sqrt n$ is large, while most of the spectral weight of $f$ is small.    
We note that if $f$ is in {\bf P} the pivotal distribution can be efficiently sampled. The spectral 
distribution can be efficiently sampled on a quantum computer (Section \ref {s:fs}).


\subsubsection {First passage percolation}

Consider an infinite planar grid where every edge is assigned a length: 1 with
probability 1/2 and 2 with probability 1/2 (independently). This model of a random metric on
the planar grid is called first-passage percolation. An old question is to understand what is the variance $V(n)$ of the
distance $D$ from $(0,0)$ to $(n,0)$? Now, let $M$ be the median value of $D$ and consider the Boolean function $f$ describing the event
``$D \ge M$''. Is $f$ noise sensitive? 

Benjamini, Kalai and Schramm (2003) showed, more or less, 
that $f$ is sensitive to logarithmic level of noise, and concluded that $V(n)=O(n/\log n)$. 
(The argument uses hypercontractivity and is similar to the argument
for critical planar percolation.) 
To show that $f$ is sensitive to noise level of $n^\delta$ for $\delta>0$
would imply that $V(n)=O(n^{1-c})$. A very interesting question is
whether methods used for critical planar percolation
for obtaining stronger noise sensitivity results can also be applied here.







\subsubsection { A beautiful problem by Benjamini and Brieussel.}

Consider $n$ steps simple random walk (SRW) $X_n$ on a Cayley graph of a finitely generated infinite group $\Gamma$. 
Refresh   independently  each step with probability $\epsilon$, to get $Y_n$ from $X_n$.
Are there groups for which the positions at time $n$, $X_n$ and $Y_n$ are asymptotically independent?
That is, the $l_1$ (total variation) distance between the 
chain  $(X_n, Y_n)$ and two independent copies $(X'_n, X''_n)$
is going to 0, with $n$.

Note that on the line $\mathbb Z$, they are uniformally correlated, 
and therefore  also on any group with a non trivial homomorphism 
to $\mathbb R$, or any group that has a finite index 
subgroup with a non trivial homomorphism to $\mathbb R$.
On the free group and for any non-Liouville group, $X_n$ and $Y_n$ are correlated as well,
but for a different reason: Both $X_n$ and $Y_n$ have nontrivial correlation with $X_1$.
Itai Benjamini  and Jeremie Brieussel 
conjecture that these are the only ways not to be noise sensitive.
That is, if a Cayley graph is Liouville and the group does not have a finite index subgroup
with a homomorphism to the reals, then the Cayley graph is noise sensitive for the simple random walk.
In particular, the Grigorchuk group is noise sensitive for the simple random walk!

\subsection {Boolean complexity, Fourier and noise}
\subsubsection { {\bf P} $\ne$ {\bf NP} --  circuit version}
The  {\bf P} $\ne$ {\bf NP}-conjecture
(in a slightly stronger form) asserts that the Boolean function described by
the graph property of containing a Hamiltonian cycle, cannot be
described by a polynomial-size circuit. 
Equivalently, the circuit form of the ${\bf NP \ne P}$-conjecture asserts that there are Boolean functions
that can be described by polynomial size nondeterministic circuits, namely as the 
projection to $n$ variables of a polynomial-size circuit, 
but cannot be 
described by polynomial size circuits.
A Boolean function $f$ is in {\bf co-NP} if $-f$ is in {\bf NP}.  

\begin {remark} Projection to $n$ variables of a Boolean function in {\bf co-NP} 
is believed to enlarge the family of functions even further. The resulting class is denoted by $\Pi_P^2$ 
and the class of functions $-f$ when 
$f \in \Pi_P^2$ is denoted by $\Sigma_P^2$. By repeating the process of negating and projecting 
we reach a whole hierarchy of complexity classes, {\bf PH}, called the polynomial hierarchy. 
\end {remark}

\subsubsection {Well below {\bf P}} 
The class {\bf NC} describes Boolean functions that 
can be expressed by polynomial size polylogarithmical depth Boolean circuits. 
This class (among others) is used to model the notion of parallel computing.
Considerably below, 
the class {\bf AC}$^0$ describes Boolean functions that can be 
expressed by bounded-depth polynomial-size circuits, where we allow AND and OR gates to 
apply on more than two inputs. A celebrated result in computational complexity asserts that 
majority and parity do not belong to {\bf AC}$^0$. 
However, the noise-stability of majority
implies that majority can be well approximated by functions in {\bf AC}$^0$.
We note that functions in {\bf AC}$^0$ are already very 
complex mathematical objects. 

A monotone threshold circuit is a circuit built from gates which are are weighted majority functions 
(without negations). 
A general threshold circuit is a circuit built from  gates which are threshold 
linear functions, i.e. we allow negative weights. {\bf TC}$^0$  ({\bf MTC}$^0$) is the class of functions 
described by bounded depth polynomial size 
(monotone) threshold circuits. 

\subsubsection {Some conjectures on noise sensitivity and bounded depth monotone threshold circuits}

\begin {conjecture}[Benjamini, Kalai, and Schramm (1999)]
\label {c:tresh}

(i) Let $f$ be a Boolean function described by a monotone threshold circuit of size $M$ 
and depth $D$. Then $f$ is stable to $(1/t)$-noise where $ t=(\log M)^{100D}$.

(ii) Let $f$ be a monotone Boolean function described by a threshold circuit of size $M$ and depth $D$. 
Then $f$ is stable to $(1/t)$-noise where $t=(\log M)^{100D}$.

\end {conjecture}

The constant 100 in the exponent is, of course, negotiable. 
In fact, replacing $ 100D$ with any function 
of $ D$ will be sufficient for most applications. 
The best we can hope for is that the conjectures are true if $t$ behaves like 
$ t=(\log M)^{D-1}$. Part (i) is plausible but looks very difficult. Part (ii) is quite reckless and may well be false. 
(See, however, Problem \ref{p:posmon}, below.) 
Note that the two parts differ ``only'' in the location of the word ``monotone.''


There are many Boolean functions that are very noise sensitive. 
A simple example is the recursive majority on threes, denoted by RM3 
and defined as follows: Suppose that $ n=3^m$. 
Divide the variables into three equal parts. Compute the RM3 separately for each of these parts and apply 
majority to the three outcomes.
Conjecture \ref {c:tresh} would have the following corollaries (C1)--(C4). 
Part (i) implies: (C1)-- RM3 is not in {\bf MTC}$^0$, and even C2 -- 
RM3 cannot be approximated by a function in Monotone {\bf MTC}$^0$.
(A variant of (C1) is known by results of Yao (1989) and Goldmann and Hastad (1991), 
and these results motivated our conjecture.) 
(C2) already seems well beyond reach.
Part (ii) implies: (C3)-- RM3 is not in  {\bf TC}$^0$ and  (C4)-- RM3 cannot 
be approximated by a function in {\bf TC}$^0$.
(We can replace RM3 with other noise-sensitive properties like 
the crossing event in planar percolation.)


\subsubsection {Bounded depth Boolean circuits and the reverse Hastad conjecture}

For a monotone Boolean function $f$ on $\Omega_n$ a Fourier description of the 
total influence is $I(f)=\sum \hat f^2(S)|S|$, and we can take this expression as the definition of $I(f)$ 
for non-monotone functions as well. The following theorem describes  briefly the situation for 
{\bf AC}$^0$. The second and third items are based on Hastad's switching lemma. 

\begin {theorem} 
(i) (Boppana (1984)): If $f$ is a  (monotone) Boolean function 
that can be described by a depth $D$ size $M$ 
monotone Boolean circuit then $ I(f) \le C(\log M)^{D-1}$.

(ii) (Hastad (1989) and Boppana (1997)) If f is a function that can be 
described by a depth $D$ size $M$ Boolean circuit then $ I(f) \le C(\log M)^{D-1}$.

(iii) (Linial Mansour Nisan (1993); improved by Hastad (2001)): If $f$ is a function that 
can be described by a depth $D$ size $M$ 
monotone Boolean circuit then $ \{\sum \hat f^2(S):|S|=t\}$ decays exponentially with $ t$ when $ t>C(\log M)^{D-1}$. 
\end {theorem}


We conjecture that functions with low influence can be approximated by low-depth small size circuits. 
A function $g$ $\epsilon$-approximates a function $f$ if $\mathbb E(f-g)^2 \le \epsilon$.

\begin {conjecture} [Benjamini, Kalai, and Schramm (1999)]
For some absolute constant $C$ the following holds. 
A Boolean function $f$ can be $0.01$-approximated 
by a circuit of depth $d$ of size $M$ where 
$(\log M)^{Cd} \le I(f).$
\end {conjecture}

\subsubsection {Positive vs. Monotone}

We stated plausible while hard conjectures on functions in  {\bf MTC}$^0$ and reckless perhaps wrong conjectures 
on monotone functions in {\bf TC}$^0$. But we cannot  present a single example of a monotone function 
in {\bf TC}$^0$ that is not in {\bf MTC}$^0$. To separate the conjectures we need 
monotone functions in {\bf TC}$^0$ that cannot even be approximated in  {\bf MTC}$^0$.
Ajtai and Gurevich (1987) proved that there are monotone functions in {\bf AC}$^0$ that are not in 
monotone {\bf AC}$^0$.


\begin {problem}
\label{p:posmon}
(i) Are there monotone functions in 
{\bf AC}$^0$ 
that cannot be approximated by functions in Monotone {\bf AC}$^0$ ?

(ii) Are there monotone functions in {\bf TC}$^0$ that are not in {\bf MTC}$^0$?

(iii) Are there monotone functions in {\bf TC}$^0$ that cannot be approximated by functions in  {\bf MTC}$^0$?
\end {problem}



\subsection { A small taste of PCP,  hardness of approximation, and \textsc {Max Cut}}
\label{s:cut}
A vertex cover of a graph $G$ is a set of vertices such that every edge contains a vertex in the set.
\textsc {Vertex Cover} 
is the algorithmic problem of finding 
such a set of vertices of minimum size. Famously this problem is an
{\bf NP}-complete problem, in fact, it is one of the problems in Karp's original list.
A matching in a graph is a set of edges such that every vertex is included in at most one edge.
Given a graph $G$ there is an easy efficient algorithm to find a maximal matching.
Finding a maximal matching with $r$ edges with respect to inclusion,
gives us at the same time a vertex cover of size $2r$ and a guarantee
that the minimum size of a vertex cover is at least $r$. 
A very natural question is to
find an efficient algorithm for a better approximation. There is by now 
good evidence that this might not be possible. It is known to 
derive (Khot and Regev (2003)) from Khot's unique game conjecture (Khot (2002)).


A cut in a graph is a partition of the vertices into two sets. The \textsc{Max Cut} problem is the
problem of finding a cut with the maximum number of edges between the parts. Also this problem 
is {\bf NP}-complete, and in Karp's list. 
The famous Goemans-Williamson algorithm
based on semidefinite programming achieves $\alpha$-approximation for max cut where 
$\alpha_ {GM} = .878567$. 
Is there an efficient algorithm for a better approximation? There is by now 
good evidence that this might not be possible.

\subsubsection {Unique games, the unique game conjecture, and the PCP theorem}
We have a connected
graph and we want to color it with $n$ colors. 
For every edge $e$ we are given an orientation of the edge and 
a permutation $\pi_e$ on 
the set of colors. In a good coloring of the edge if the tail is colored $c$ then the head must be colored $\pi_e (c)$.
It is easy to check efficiently if a global good coloring exists since coloring one vertex forces the coloring of all others.

Given $\epsilon,\delta $ the unique game problem is to algorithmically decide between two scenarios 
(when we are promised that one of them holds.) Given a graph, $G$, a color set of size $n$ and 
a permutation constraint for each edge.

(i)  There is no coloring with more than $\epsilon$ fraction of the edges are colored good.  

(ii)  There is a coloring for which at least fraction $1-\delta$ of the edges are colored good.  

The unique game conjecture asserts that for every $\epsilon >0$ and $\delta >0$  
it is NP-hard to decide between these two scenarios.

If one does not insist on the constraints being permutations and instead 
allows them to be of general form, then the above holds, and is called the PCP Theorem --
one of the most celebrated theorems in 
the theory of computation. 


\begin {remark} A useful way to describe the situation (which also reflects the 
historical path leading to it) is 
in terms of a three-player game - there are two ``provers'' and a verifier. A verifier is 
trying to decide which of the two cases he is in, and can communicate 
with two all powerful (non-communicating) provers. To do that, the verifier 
samples an edge, and sends one endpoint to each prover. Upon receiving their answers, 
the verifier checks that the two colors satisfy the constraint.  
The provers need to convince the verifier that 
a coloring exists by giving consistent answers to simultanous questions drawn at random.  
\end {remark}

\subsubsection {The theorem of Khot, Kindler, Mossel, and O'Donell}

\begin {theorem}
\label {t:kkmo}
[Khot, Kindler, Mossel, O'Donnell (2007)]
Let $\beta > \alpha_{GM}$ 
be a constant. Then an efficient $\beta$-approximation algorithm for 
\textsc{Max Cut}  implies an efficient algorithm for unique-games.
\end {theorem}

The reduction relies on the majority is stablest theorem (Theorem \ref {t:moo}) 
which was posed by Khot, Kindler, Mossel, and O'Donnell as a conjecture and later proved by 
Mossel, O'Donnell, and Oleszkiewicz (Theorem \ref{t:moo}). 
This result belongs to the theory of hardness of approximation and probabilistically 
checkable proofs which is among the 
most important areas developed in computational complexity in the last three decades. 
For quite a few problems in Karp's original 
list of NP-complete problems (and many other problems added to the list), 
there is good evidence that the best efficient approximation 
is achieved by a known relatively simple algorithm. For a large class of problems it is even known (Raghavendra (2008)) 
(based on hardness of the unique game problem) that the best algorithm is either a very simple combinatorial algorithm 
(like for \textsc{Vertex Cover}), 
or a  more sophisticated application of semidefinite programming (like for \textsc {Max Cut}).
I will give a quick and very fragmented taste on three ingredients of the proof of Theorem \ref {t:kkmo}.

{\it The noisy graph of the cube}
The proof of the hardness of max cut relative to unique games is based on the weighted graph 
whose vertices are the vertices of the discrete cube, all pairs are edges, and the weight of 
an edge between two vertices of distance $k$ is $(1-p)^kp^{n-k}$. It turns out that in order to 
analyze the reduction, it suffices to study the structure of good cuts in this very special graph.

{\it The least efficient error correcting codes}.
Error correcting codes have, for many decades, been among the most celebrated 
applications of mathematics with huge impact on technology. 
They also have a prominent role in 
theoretical computer science and in PCP theory.
The particular code needed for max cut is the following:
Encode a number $k$ between 1 to $n$ (and thus $\log n$ bits) by a Boolean function - 
a dictatorship where the $k$th variable is the dictator!

{\it Testing dictatorship}
An important ingredient of a PCP proof is ``property testing'', testing by 
looking at a bounded number of values if a Boolean function satisfies 
a certain property, or is very far from satisfying it. In our case we would like to test 
(with high probability of success) if a Boolean function  
is very far from dictatorship, 
or has substantial correlation 
with it. The test is the following:  
Choose $x$ at random, let $y=N_\epsilon (-x)$. Test if $f(x) = -f(y)$.
For the majority function 
the probability that  majority passes the test is roughly $\arccos(\epsilon -1)$, 
majority is stablest theorem  implies that anything that is more stable 
has large correlation with a dictator.

\subsubsection {Discussion: integrality gap and polytope integrality gap}

Given a graph $G$ and nonnegative weights on its vertices, 
the weighted version of vertex cover is the algorithmic problem of finding 
a set of vertices of minimum weight that covers all edges. 

{\bf Minimize}  $w_1x_1+w_2x_2+\cdots +w_nx_n$ where $x=(x_1,x_2,\dots,x_n)$ is a 0-1 vectors, 

{\bf subject to:} $x_i+x_j \ge 1$ for every edge $\{i,j\}$.

Of course, this more general problem is also {\bf NP}-complete. 
The linear programming relaxation allows $x_i$s to be real  
belonging to the interval [0,1].
The {\it integrality gap} for general vertex cover problems is 2 and given the solution to 
the linear programming problem you can just consider 
the set of vertices $i$ so that $x_i \ge 1/2$. 
This will be a cover and the ratio between this cover and the optimal one is at most 2.
The integrality gap for the standard relaxation of max cut is $\log n$. The integrality gap is 
an important part of the picture in PCP theory. I conclude with a beautiful problem that 
I learned from Anna Karlin.
 
Consider the integrality gap (called the {\it polytope integrality gap}) 
between the covering problem and the linear programming relaxation 
when the graph $G$ is fixed. In greater generality, consider a general 
covering problem of maximizing $c^tx$ subject to $Ax \le b$ where $A$ 
is integral matrix  of nonnegative integers. 
Next,  considered the integrality gap between 0-1 solutions and real solutions 
in $[0,1]$ when $A$ and $b$ are fixed (thus the feasible polyhedron is fixed, hence 
the name ``polytope integrality gap'') 
and only $c$ (the objective function) varies. 
The problem is  
if for vertex cover for every graph $G$ and every vector of weights,  
there is an efficient algorithm achieving the polytope integrality gap. 
The same question can be asked
for polytope integrality gap of arbitrary covering problems.

\section {The quantum computer challenge}
\label {s:qc}

\begin {center} 
{\color {blue} 
{\it To Robert Aumann, Maya Bar-Hillel, Dror Bar-Nathan, Brendan McKay and
Ilya Rips who trained me as an applied mathematician.}
}
\end {center}

\subsection {Quantum computers and noise}
Recall that the basic memory component in classical computing is a ``bit,'' which
can be in two states, ``0'' or ``1.'' A computer, as modeled by a Boolean circuit, has
$n$ bits and it can perform certain logical operations  on them. 
The NOT gate, acting on a single bit, and the AND
gate, acting on two bits, suffice for {\it universal} classical computing.
This means that a computation based on another collection of logical gates,
each acting on a bounded number of bits,
can be replaced by a computation based only on NOT and AND.
Classical circuits equipped with random bits lead to randomized algorithms, which, as mentioned before, 
are both practically useful
and theoretically important.
Quantum computers allow the creation of probability
distributions that are well beyond the reach of classical computers with access to random bits.

\subsubsection {Quantum circuits}
A qubit is a piece of quantum memory.
The state of a qubit can be described by a unit
vector in a two-dimensional complex Hilbert space $H$.
For example, a basis for $H$ 
can correspond to two energy levels of the hydrogen
atom, or to horizontal and vertical polarizations of a photon.
Quantum mechanics allows the qubit to be in a {\it superposition} of the basis vectors, described by
an arbitrary unit vector in $H$.
The memory of a
quantum computer (``quantum circuit'') consists of $n$ qubits. Let $H_k$ be the two-dimensional Hilbert space
associated with the
$k$th qubit.
The state of the entire memory of $n$ qubits is described by
a unit vector in the tensor product $H_1 \otimes H_2 \otimes \cdots \otimes H_n$.
We can put one
or two
qubits through {\it gates} representing
unitary transformations acting on
the corresponding two- or four-dimensional Hilbert spaces, and as for classical computers, there is a
small list of gates sufficient for  universal quantum computing.
At the end of the computation process, the state of the entire computer can be {\it measured}, giving a probability
distribution on 0--1 vectors of length $n$.


A few words on the connection between the mathematical model of quantum circuits and quantum physics:
in quantum physics, states and
their evolutions (the way they change in time)
are governed by the Schr\"odinger equation.
A solution of the
Schr\"odinger equation can be described
as a unitary process on a Hilbert space and
quantum computing processes of the kind we just described form a large class of such quantum evolutions.

\begin {remark} Several universal classes of quantum gates are described in 
Nielsen and Chuang (2000) [Ch. 4.5].
The gates for the IBM quantum computer are eight very basic one-qubit gates, 
and the 2-qubit CNOT gate 
according to a certain fixed directed graph. This is a universal system 
and in fact, an over complete one.  
\end {remark}


\subsubsection {Noise and fault-tolerant computation}
\label {s:noise}

The main concern
regarding the feasibility of quantum computers has always been that
quantum systems are
inherently noisy: we cannot accurately control them, and we
cannot accurately describe them. The concern regarding noise in quantum systems as a major
obstacle to quantum computers was put forward
in the mid-90s by Landauer (1995), Unruh (1995), and others.

What is noise? As we said already, solutions of the Schr\"odinger equation (``quantum evolutions'') 
can be regarded as unitary
processes on Hilbert spaces.
Mathematically speaking, the study of noisy quantum systems is the
study of {\it pairs} of Hilbert spaces $(H, H')$, $H \subset H'$,
and a unitary process on the larger Hilbert space $H'$.
Noise refers to the general effect of neglecting degrees
of freedom, namely, approximating the process on a large Hilbert space
by a process on the small Hilbert space. For controlled quantum systems and,
in particular, quantum computers, $H$ represents the controlled part of the system,  and
the large unitary process on $H'$  represents, in addition to an ``intended''
controlled evolution on $H$,
also the uncontrolled effects of the {\em environment}.
The study of noise is relevant, not only to controlled
quantum systems, but also to many other aspects of quantum physics.

A second, mathematically equivalent way
to view noisy states and noisy evolutions, is to stay with the original Hilbert space $H$, but to
consider a mathematically larger class of states and operations.
In this view, the state of a noisy qubit is described as a classical probability distribution on
unit vectors of the associated Hilbert spaces. Such states are referred to as {\it mixed states}.

It is convenient to think about the following simple form of noise, called {\it depolarizing noise}:
in every computer cycle a qubit is not affected with probability $1-p$, and, with probability
$p$, it turns into the \emph{maximal entropy mixed state}, i.e.,
the average of all unit vectors in the associated Hilbert
space.

\begin {remark} 
It is useful to distinguish between the model error rate which is $p$  
in the above example, 
and the effective error rate: the 
probability that 
a qubit is corrupted
at a computation step, conditioned on it having survived up to this step.
The effective error rate depends not only on the model error rate but 
also on the computation sequence. When the computation is non trivial (for example,
for pseudo random circuits) the effective error rate grows linearly with the 
number of qubits. This is a familiar fact that is taken into account by the threshold 
theorem described below. 
\end {remark}

To overcome noise, a 
theory of quantum fault-tolerant computation based on quantum error-correcting codes
was developed. 
Fault-tolerant computation refers to computation in the presence of errors. The
basic idea is to represent (or ``encode'') a single (logical) qubit
by a large number of physical qubits, so as to ensure
that the computation is robust even if some of these physical qubits are faulty.

\begin {theorem} 
[Threshold theorem (informal statement) (Aharonov, Ben-Or (1997), Kitaev (1997), Knill, Laflamme, Zurek (1998)]
\label {t:tt} 
When the level of noise is below a certain positive threshold $\rho$, 
noisy quantum computers allow universal quantum computation.
\end {theorem}


Theorem \ref {t:tt} shows that once high quality quantum circuits are  built for roughly 100-500 qubits 
then it will be possible in principle to use quantum error-correction codes 
to amplify this achievement for building quantum computers 
with unlimited  number of qubits. The interpretation of this result took for granted that quantum computers with a 
few tens of qubits are feasible, and this is incorrect.

Let $A$ be the maximal number of qubits for which a reliable quantum circuit can be engineered.
Let $B$ be the number of qubits required for good quantum error-correcting codes needed 
for quantum fault tolerance. $B$ is in the range of 100-1000 qubits.

{\bf The optimistic scenario} $A > B$ 

{\bf The pessimistic scenario} $B > A$

As we will see, there are good theoretical reasons for 
the pessimistic scenario (even for $B \gg A$) as well as interesting consequences from it.
We emphasize that both scenarios are compatible with quantum mechanics.


\begin{figure}
\centering
\includegraphics[scale=0.4]{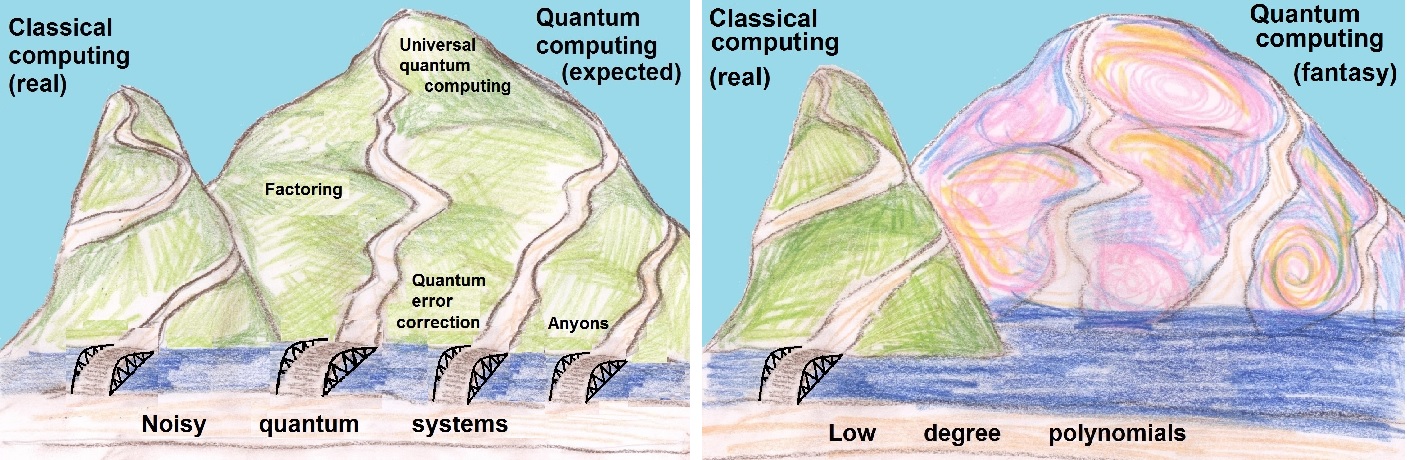}
\caption{\Small {{\bf Two scenarios:}  {\bf Left} --  
Quantum fault-tolerance mechanisms, via quantum error-correction, allow robust quantum information and
computationally superior quantum computation.
{\bf Right} -- Noisy quantum evolutions,
described by low-degree polynomials, allow, via the
mechanisms of averaging/repetition,
robust {\it classical} information and computation,  but do not allow reaching the
starting points for quantum supremacy and quantum fault-tolerance.
 Drawing by Neta Kalai.}}
\label{fig:OP}
\end{figure}




\subsection {Complexity 6: quantum computational supremacy}



\subsubsection {\textsc {Factoring} is in {\bf BQP}}

Recall that computational complexity
is the theory of {\it efficient computations}, where ``efficient''
is an asymptotic notion referring to situations where
the number of computation steps (``time'') is at most a polynomial in the number of input bits.
We already discussed the complexity classes {\bf P} and {\bf NP}, and let us (abuse notation and) allow 
to add classical randomization to all complexity classes we discuss.
There are important intermediate problems between {\bf P} and {\bf NP}. \textsc{Factoring} -- the task of 
factoring an $n$-digit integer to its prime decomposition
is not known to be in {\bf P}, as the best algorithms
are exponential in the cube root of the number of digits. \textsc{Factoring} is in {\bf NP}, hence it is unlikely
that factoring is {\bf NP}-complete. Practically, \textsc{Factoring} is hard and this is the 
basis to most of current cryptosystems. 


The class of decision
problems 
that quantum computers can efficiently solve is denoted by {\bf BQP}.
Shor's algorithm shows that quantum computers can factor $n$-digit integers efficiently --
in $\sim n^2$ steps! Quantum computers are not known to be able to solve
{\bf NP}-complete problems efficiently, and there are good reasons to think that they cannot do so. However,
quantum computers can
efficiently perform certain computational tasks
beyond {\bf NP} and even beyond {\bf PH}. 

The paper by Shor (1994) presenting an  efficient factoring algorithm for quantum computers 
is among the scientific highlights of the 
20th century, with immense impact on 
several theoretical and experimental areas of physics.

\subsubsection {Fourier sampling, boson sampling and other quantum sampling} 

\label {s:fs}

As mentioned in Section \ref {s:lp} exact 
and approximate sampling are important algorithmic tasks on their own and as subroutines for other tasks. 
Quantum computers enable remarkable new forms of sampling. 
Quantum computers would allow the creation of probability
distributions that are beyond the reach of classical computers with access to random bits.
Let \textsc{Quantum Sampling} 
denote the class of distributions that quantum computers can efficiently sample.
An important class of such distributions 
is \textsc{Fourier Sampling}.
Start with a Boolean function $f$. (We can think about $f(x)$ as the winner in HEX.) 
If $f$ is in {\bf P} then 
we can classically sample $f(x)$ for a random $x \in \Omega_n$. With quantum computers we can do more. 
A crucial ability of quantum computers is to prepare a state
 $2^{-n/2} \cdot \sum f(x)|x>$ which is a superposition of all $2^n$ vectors weighted by the value of $f$. 
Next a quantum computer can easily take the Fourier transform 
of $f$ and thus sample exactly a subset $S$ according to $\hat f^2(S)$.
This ability of quantum computers goes back essentially to Simon (1995), and is crucial for Shor's factoring algorithm.

Another important example is  \textsc{Boson Sampling} that refers to 
a  class of probability distributions representing a collection of non-interacting bosons,
that quantum computers can efficiently create. 
\textsc{Boson Sampling} was introduced by Troyansky and Tishby (1996) and was intensively studied
by Aaronson and Arkhipov (2013), 
who offered it as a quick path for experimentally showing that quantum supremacy is a real phenomenon.
Given an $n$ by $n $ matrix $A$, let  $per (A)$
denote the permanent of $A$. 
Let $M$ be a complex $n \times m$ matrix, $m \ge n$ with orthonormal rows.
Consider all ${m+n-1} \choose {n}$ sub-multisets $S$ of $n$ columns (namely, allow columns to repeat),
and for every sub-multiset $S$ consider the corresponding $n \times n$ submatrix $A$ (with column $i$ repeating $r_i$ times).
\textsc {Boson Sampling} 
is the algorithmic task of sampling those multisets $S$ according to
$|per (A)|^2/(r_1!r_2!\cdots r_n!)$. 

\subsubsection {Hierarchy collapse theorems} 

Starting with Terhal and DiVincenzo (2004) there is a series of works showing that it is very unreasonable to expect that 
a classical computer can perform quantum sampling even regarding distributions that express very limited quantum computing.
(Of course quantum computers can perform these sampling tasks).

\begin {theorem} [Terhal and DiVincenzo (2004), Aaronson and Arkhipov (2013), 
Bremner, Jozsa, and Shepherd (2011), and others ] 
\label{t:qs}

If a classical computer can exactly sample according to either

(i) General \textsc{Quantum Sampling},

(ii) \textsc{Boson Sampling},

(iii) \textsc{Fourier Sampling},

(iv) Probability distributions obtained by bounded 
depth polynomial size quantum circuits --

then the polynomial hierarchy collapses. 
\end {theorem}

The basic argument is to show that if a classical computer 
which allows any of the sampling tasks listed above, is equipped with an {\bf NP} oracle, 
then  it is able to efficiently perform {\bf \#P-complete} computations.
This will show that the class {\bf \#P} that includes {\bf PH} already collapses
to the third level in the polynomial hierarchy.



\subsection {Computation and physics 1}


 
\subsubsection{Variants on the Church-Turing thesis}

The famous Church-Turing thesis (CTT) asserts that everything computable is computable by a Turing machine. 
Although initially this was a thesis about computability, there were early attempts to relate it to physics, 
namely to assert that physical devices obey the CTT. 
The efficient (or strong) Church-Turing 
thesis (ECTT) in the context of feasible computations by physical devices was considered early on by 
Wolfram (1985), Pitowski (1990) and others. It asserts that only efficient computations by a Turing machine 
are feasible physical computation. Quantum computers violate the ECTT. 
The pessimistic scenario brings us back to the ECTT, 
and, in addition, it proposes 
even stronger limitation for ``purely quantum processes'' (suggested from Kalai, Kindler (2014)).

\begin {quotation}
Unitary evolutions that can be well-approximated by physical devices, can be approximated by low degree polynomials, 
and are efficiently learnable.
\end {quotation}
 
The following {\bf NPBS}-principle (no primitive-based supremacy) 
seems largely applicable in the interface between practice and theory 
in the theory of computing.

\begin {quotation} 
Devices that express (asymptotically) primitive (low-level) computational power 
cannot be engineered or programmed to achieve superior computational tasks.
\end {quotation}

\subsubsection {What even quantum computers cannot achieve and the modeling of locality}

The model of quantum computers already suggests important limitations on what local quantum systems can compute. 

\begin {itemize}
\item
{\it Random unitary operations on large Hilbert spaces.}  
A quantum computer with $n$ qubits cannot reach a random unitary state since reaching such a state requires an 
exponential number of 
computer cycles. (Note also that since an $\epsilon$-net of states for $n$-qubits quantum computer requires a set of 
size doubly exponential in $n$, most states are beyond reach for a quantum computer.)  
\item
{\it Reaching ground states for complex quantum systems.}
A quantum computer is unlikely to be able to reach the ground 
state of a quantum system (that admit an efficient description).
As a matter of fact, reaching a ground state is {\bf NP}-complete 
even for classical systems and for quantum computing
the relevant complexity class is even larger {\bf QMA}.
\end {itemize}


These limitations are based on the model of quantum computers (and the second also on 
{\bf NP $\ne$ P}) 
and thus do not formally 
follow from the basic framework of quantum mechanics (for all we know). They do follow 
from a principle of ``locality'' asserting that quantum evolutions express interactions between a 
small number of 
physical elements. This principle is modeled by quantum computers, 
and indeed a crucial issue in the debate on quantum computers
is what is the correct modeling of local quantum systems. 
Let me mention three possibilities.

\begin {enumerate}
\item[(A)]

The model of quantum circuits is the correct model for local quantum evolutions. 
Quantum computers are possible, the difficulties are matters of engineering,  and 
quantum computational supremacy is amply manifested in quantum physics.

\item [(B)]
The model of noisy quantum circuits is the correct model for local quantum evolutions. 
In view of the threshold theorem, quantum computers are possible 
and the remaining difficulties are matters of engineering.

\item [(C)]
The model of noisy quantum circuits is the correct model for local quantum evolutions, and 
further analysis suggests that the threshold in the threshold theorem cannot be reached.  
Quantum circuits with noise above the threshold is the correct modeling of local quantum systems.  
Quantum computational supremacy is an artifact of incorrect modeling of locality.
\end {enumerate}

Computational complexity insights (and some common sense) can assist us deciding between these possibilities. 
While each of them 
has its own difficulties, in my view the third one is correct.

\subsubsection {Feynman's motivation for quantum computing}

\begin {conjecture}[Feynman's (1982) motivation for quantum computation]
High energy physics computations, 
especially computations in QED (quantum electrodynamics) 
and QCD (quantum chromodynamics), can be carried out efficiently by quantum computers.
\end {conjecture}

This question touches on the important mathematical question of giving rigorous mathematical foundations 
for QED and QCD 
computations. 
Efficient quantum computation for them 
will be an important (while indirect) step 
toward putting these theories on rigorous mathematical grounds. 
Jordan, Lee, and Preskill (2012, 2014) found an efficient algorithm for certain 
computations in ($\phi^4$) quantum field theory for
cases where a rigorous mathematical framework is available.


\subsection {The low scale analysis: Why quantum computers cannot work}

\subsubsection {Noisy systems of non interacting photons} 

\begin{figure}
\centering
\includegraphics[scale=0.6]{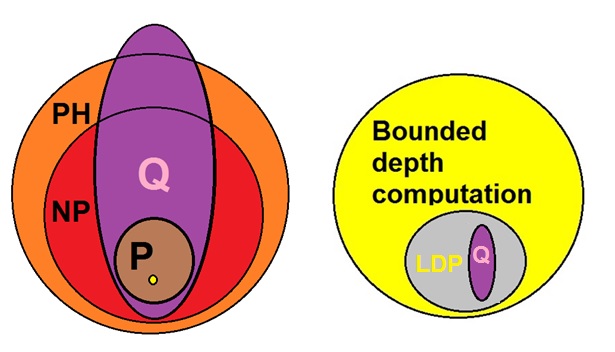}
\caption{{\Small  
Quantum computers offer mind-boggling computational superiority (left), but in the small 
scale, noisy quantum circuits are computationally very weak, 
unlikely to allow quantum codes needed for quantum computers (right). 
}}
\label{fig:3.2}
\end{figure}

\begin{theorem}[Kalai and Kindler (2014)]
\label {thm:1}
When the noise level is constant, \textsc{Boson Sampling} distributions are well approximated
by their low-degree Fourier--Hermite expansion.
Consequently, noisy \textsc {Boson Sampling} can be approximated by
bounded-depth polynomial-size circuits.
\end {theorem}

It is reasonable to assume that for all proposed implementations of \textsc {Boson Sampling}
the noise level is at least a constant, and therefore, an experimental realization
of \textsc {Boson Sampling} represents, asymptotically, bounded-depth computation. 
In fact 
noisy \textsc {Boson Sampling} belongs to a computational class {\bf LDP}  (approximately sampling distributions  
described by bounded degree polynomials) which is well below {\bf AC}$^0$.
The next theorem shows that
implementation of \textsc {Boson Sampling} will actually require pushing down the noise level
to below $1/n$.

\begin{theorem}[Kalai and Kindler (2014)]
\label {thm:2}
When the noise level is $\omega(1/n)$, and $m \gg n^2$,
\textsc {Boson Sampling}  is very sensitive to noise with a
vanishing correlation between the noisy distribution
and the ideal distribution.\footnote{The condition $m \gg n^2$ can
probably be removed by a more detailed analysis.}
\end {theorem}




\subsubsection {Noisy quantum circuits}

\begin {conjecture}
\label {c:pc}

(i) The insights for noisy \textsc{Boson Sampling} apply to all versions of 
realistic forms of noise for non interactive bosons.

(ii) These insights extend further to quantum circuits and other quantum devices in the small scale.

(iii) These insights extend even further to quantum devices, 
including microscopic processes, that do not use quantum error-correction

(iv) This deems quantum computational supremacy and the needed quantum error correction codes impossible
\end {conjecture}

The first item seems quite a reasonable extension of Theorem \ref {thm:1}. In fact, the argument 
applies with small changes to a physical modeling of mode-mismatch noise. 
(When bosons are not fully indistinguishable.) 
Each item represents  quite a leap from the previous one. The last item expresses the idea that superior computation
cannot be manifested by primitive asymptotic computational power. Theorem \ref {thm:1} 
put noisy \textsc{Boson Sampling} in a very low-level class, {\bf LDP} even well below {\bf AC}$^0$.  
It is not logically impossible but still quite implausible that such a primitive computing device will 
manifest superior computing power for 50 bosons. 

\begin {remark}
{\it Why robust classical information and computation is possible and ubiquitous.} The ability to approximate
low-degree polynomials
still supports robust {\it classical} information. This is related to our second puzzle.
The majority function
allows for very robust bits
based on a large number of noisy bits and
admits
excellent low-degree approximations. Both encoding (by some repetition procedure) 
 and decoding (by majority or a variation of majority) needed for robust 
classical information are supported by 
low degree polynomials.
\end {remark}


\subsection {Predictions regarding intermediate goals and near-term experiments}


\begin {itemize}

\item

{\it Demonstrating quantum supremacy.}
A demonstration of quantum computing supremacy, namely crossing the line where classical simulation is 
possible, requires, e.g., building of pseudo-random quantum circuits of 50-70 qubits. 
As we mentioned in the Introduction, this idea can be partially tested already for quantum circuits with 10-30 qubits, 
and there are plans for a decisive  
demonstration on 50 qubits in the near future. 
Quantum supremacy could be demonstrated 
via implementation of \textsc { Boson Sampling} and in various other ways. 
Theorems \ref{thm:1} and \ref {thm:2} and principle {\bf NPBS} suggest that all these attempts will fail. 

\item
{\it Robust quantum qubits via quantum error-correction.}
The central goal towards quantum computers is to build logical qubits 
based on quantum error-correction, and a major 
effort is made to demonstrate a distance-5 
surface code which requires 100 or so qubits. It is now commonly agreed that this task is 
harder than ``simply'' demonstrating quantum computational supremacy. Therefore, principle {\bf NPBS}
suggests that these attempts will fail as well.


\item
{\it Good quality individual qubits and gates (and anyonic qubits).}
The quality of individual qubits and gates is the major factor in the quality of quantum circuits built from them. 
The quantum computing analogue of Moore's law, known as  ``Schoelkopf's law'' asserts that 
roughly every three years, quantum decoherence can be delayed by a factor of ten.
The analysis leading to the first two items suggests 
that Schoelkopf's law will break before reaching the quality needed for quantum supremacy and
quantum fault tolerance. This is an indirect argument, but more directly, 
the microscopic process leading to the qubits also (for all we know) represents low level complexity power.
This last argument also sheds doubt on hopes of reaching robust quantum qubits via anyons.

\end {itemize}


\subsection {Computation and physics 2: 
noisy quantum systems above the noise threshold}

Basic premises for studying noisy quantum evolutions under the pessimistic scenario are first that  
the emerged modeling is implicit;  namely, it is given in terms of
conditions that the noisy process must satisfy, rather than a direct description for the noise.
Second, there are systematic relations between the (effective) noise and the entire quantum
evolution and also between the target state
and the noise.

\subsubsection {Correlated errors and error synchronization} 


The following prediction regarding noisy entangled pairs of qubits (or ``noisy cats'')
is perhaps the simplest prediction on noisy quantum circuits under the pessimistic scenario.
Entanglement is a name for quantum correlation,
and it is an important feature of
quantum physics and a crucial ingredient of quantum computation.
A \emph{cat state} of the
form ${\frac{1}{\sqrt 2}}\left|00\right\rangle +{\frac{1}{\sqrt 2}} \left|11\right\rangle$
represents the simplest (and strongest) form of entanglement between two qubits.

{\bf Prediction 1: Two-qubits behavior.}
For any implementation of quantum circuits, cat states are subject to qubit errors 
with substantial positive correlation.


Error synchronization refers to a substantial probability that a large number of qubits, 
much beyond the average rate of noise,  
are corrupted. This is a very rare phenomenon for the model noise and when quantum 
fault-tolerance is in place, 
error-synchronization is an extremely rare event also for the 
effective noise.

{\bf Prediction 2: Error synchronization.} For pseudo random circuits 
highly synchronized errors will necessarily occur.

\begin {remark}
Both predictions 1 and 2 can already be tested 
via the quantum computers of Google, IBM and others. 
(It will be interesting to test prediction 1  
even on gated qubits, where it is not in conflict with the threshold theorem, 
but may still be relevant to the 
required threshold constant.)
\end {remark}



\subsubsection {Modeling general noisy quantum systems}
\hspace{0.5cm}

{\bf Prediction 3: Bounded-degree approximations, and effective learnability.}
Unitary evolutions that can be approximated by noisy quantum circuits (and other devices)
are approximated by low degree polynomials and are efficiently learnable.  

{\bf Prediction 4: Rate.}  For a noisy quantum system a lower bound for the rate of noise in a time interval
is a measure of non-commutativity for the projections in the algebra of unitary operators
in that interval.

{\bf Prediction 5:  Convoluted time smoothing.}
Quantum evolutions are subject to noise with a substantial correlation with time-smoothed evolutions.

Time-smoothed evolutions
form an interesting restricted class of noisy quantum evolutions
aimed for modeling evolutions under the pessimistic scenario when fault-tolerance is
unavailable to suppress noise propagation.
The basic example for time-smoothing is the following:
start with an ideal quantum evolution
given by a sequence of $T$ unitary operators,
where $ U_t$ denotes the unitary
operator for the $t$-th step, $t=1,2,\dots T$.
For $s<t$ we denote $ U_{s,t} = \prod_{i=s}^{t-1}U_i$ and let $U_{s,s}=I$ and $U_{t,s}=U^{-1}_{s,t}.$
The next step is to add noise in a completely standard way: consider a noise operation $ E_t$ for
the $t$-th step. We can think about the case where the unitary evolution is a quantum computing process
and $E_t$ represents a depolarizing noise with a fixed
rate acting independently on the qubits.
And finally, replace $E_t$ with a new noise operation $E'_t$ defined as the average

\begin {equation}
\label {e:smooth}
 E'_t =
\frac{1}{T}
\cdot
\sum_{s=1}^T  U_{s,t} E_s U^{-1}_{s,t}.
\end {equation}

Predictions 1-5 are implicit and describe systematic relations between
the (effective) noise and the evolution.
We expect that time-smoothing will suppress
high terms for some Fourier-like expansion\footnote {Pauli expansion seems appropriate for the 
case of quantum circuits, see Montanaro, Osborne (2010)}, thus relating Predictions 3 and 5.
Prediction 4 resembles the picture drawn by Polterovich (2007) 
of the ``unsharpness principle'' in symplectic geometry, quantization and quantum noise. 

\begin {remark} It is reasonable to assume that time-dependent quantum evolutions are inherently noisy since 
time dependency indicates interaction with an environment. Two caveats: famously, time dependent evolutions 
can be simulated by time independent evolutions, but we can further assume that in such cases 
the noise lower bounds will transfer. Second, in the context of noise, environment of an 
electron (say) refers also to 
its internal structure.
\end {remark}

\begin {remark} 
Physical processes are not close to unitary evolutions and there are 
systematic classical effects (namely robust effects of 
interactions with a large environment.) We certainly cannot 
 model everything with low degree polynomials. On the other hand, it is unlikely that natural 
physical evolutions express the full power of {\bf P}. It will be interesting to understand the complexity of various 
realistic physical 
evolutions , and to identify larger relevant classes within {\bf P}, especially classes for which efficient 
learnability is possible. (Efficient learnability is likely to be possible in wide contexts. 
After all, we have learned quite a bit.) 
\end{remark}



\begin {remark}
Let me list (for more details see Kalai (2016, expanded)) 
a few features of noisy quantum evolutions and states ``above the threshold.''
(a) {\it Symmetry.} Noisy quantum states and evolutions are subject to noise that respects
their symmetries. (b) {\it Entropy lower bounds.} Within a symmetry class of quantum states/evolutions (or
for classes of states defined in a different way), there is an absolute positive lower bound for entropy.
(c) {\it Geometry.} Quantum states and evolutions reveal
some information on the geometry of (all) their physical realizations.
(d) {\it Fluctuation.}
Fluctuations in the rate of noise for interacting $N$-element systems
(even in cases where interactions
are weak and unintended) scale like $N$ and not like $\sqrt N$.
(e) {\it Time.} The difficulty in implementing a local quantum computing process is not invariant under reversing time.
\end {remark}


\subsection {Noster computationalis mundus (Our computational world)}
The emerging picture from our analysis is that the 
basic computational power of quantum devices is very limited: 
unitary evolutions described by noisy local quantum devices  
are confined to low degree polynomials. It is 
classical information and computation which emerge via noise-stable encoding and decoding processes 
that enable the wealth of computation witnessed in nature. 
This picture offers many challenges, in the mathematical, physical, and computational aspects, and
those can serve as a poor man's replacement for quantum supremacy dreams. 

\section {Conclusion}

We talked about three fascinating puzzles on mathematics and
computation, telling a story which involves pure
and applied mathematics, theoretical computer science, 
games of various kinds, physics, and social sciences. The connection and tension between
the pure and the applied, between models and reality, and 
the wide spectrum between foundations and engineering
is common to all of our three puzzles.
We find great expectations, surprises, mistakes, disappointments
and controversies
 at the heart of our endeavor, while seeking truth and understanding
in our logical, physical,
and human reality. 
In our  sweet professional lives, 
being wrong 
while pursuing dreams unfounded in reality, is sometimes of value, 
and second only to being right. 



\section* {Appendix: abstract objective functions and
telling a polytope from its graph} 

The linear programming local-to-global principle has very nice connections and applications
to the combinatorial theory of convex polytopes. Abstract linear objective
functions (and Sharir-Welzl's abstract linear programming problems)
are related to the notion of shellability. We will bring here one such application Kalai (1988) -- 
a beautiful proof
that I found to the following theorem of 
Blind and Mani (1987) conjectured by Micha A. Perles.

\begin {theorem}
\label {t:bm}
The combinatorial structure of a simple polytope $P$ is determined by its graph.
\end{theorem}

We recal that a $d$-polytope $P$ is simple if every vertex belongs to exactly $d$ edges.
Thus the graph of $P$ is a $d$-regular graph. For a simple polytope, 
every set of $r$ edges containing a vertex $v$ 
determines an $r$-face of $P$. Faces of simple polytopes are simple.
Consider an ordering $\prec$ of the vertices of a simple $d$-polytope $P$.
For a nonempty face $F$ we say
that a vertex $v$ of $F$ is a local maximum in $F$ if $v$ is larger w.r.t. the
ordering $\prec$ than all its neighboring vertices in $F$.
Recall that  an abstract objective function (AOF) of a simple $d$-polytope
is an ordering which satisfies the basic property
of linear objective functions:
Every nonempty face $F$ of $P$ has a unique local maximum vertex.

If $P$ is a simple $d$-polytope and $\prec$ is a linear ordering of the vertices
we define the degree of a vertex $v$ w.r.t. the ordering
as the number of adjacent vertices to $v$ that are smaller than $v$
w.r.t $\prec$. Thus, the degree of a vertex is a nonnegative number
between 0 and $d$. Let  $h_k^\prec$ be the number of vertices of degree $k$.
Finally, put $F(P)$ to be the total number of nonempty faces of $P$.

{\bf Claim 1:} $\sum_{r=0}^d 2^k h_k^\prec ~~~\ge~~~ F(P),$
and equality hold if and only if
the ordering $\prec$ is a AOF.

{\bf Proof:} Count pairs $(F,v)$ where
$F$ is a non-empty face of $P$ (of any dimension) and $v$ is a
vertex which is local maximum in $F$ w.r.t. the ordering $\prec$.
On the one hand every vertex $v$ of degree $k$ contributes
precisely $2^k$ pairs $(F,v)$ corresponding to all subsets of edges
from $v$ leading to smaller vertices w.r.t. $\prec$. Therefore the
number of pairs is precisely $\sum_{r=0}^d 2^k h_k^\prec$.
On the other hand,
the number of such pairs is at least $F(P)$ (every face has at least
one local maximum) and it is equal to $F(P)$ iff every face has
exactly one local maximum i.e if the ordering
is an AOF.

{\bf Claim 2:} A connected $k$-regular subgraph $H$ of $G(P)$ is the graph of
a $k$-face, if and only if there is an AOF in whic
h all vertices in $H$
are smaller than all vertices not in $H$.

{\bf Proof:}
If $H$ is the graph of a $k$-face $F$ of $P$ then consider a linear objective
function $\psi$ which attains its minimum precisely at the points in $F$.
(By definition for every non-trivial face such a linear objective function
exists.) Now perturb $\psi$ a little to get a generic linear objective
function $\phi$ in which all vertices of $H$ have smaller values
than all other vertices.

On the other hand if there is an AOF $\prec$ in which all vertices in $H$
are smaller than all vertices not in $H$, consider the vertex $v$ of $H$
which is the largest w.r.t. $\prec$. There is a $k$-face $F$ of $P$ determined
by the $k$-edges in $H$ adjacent to $v$ and $v$ is a local maximum
in this face. Since the ordering is an AOF $v$ must be larger than
all vertices of $F$ hence the vertices of $F$ are contained in $H$ and
the graph of $F$ is a subgraph of $H$. But the only $k$-regular
subgraph of a connected $k$-regular graph is the graph itself and
therefore $H$ is the graph of $F$.

{\bf Proof of Theorem \ref {t:bm}:}
Claim 1 allows us to determine just from the graph all the
orderings which are AOF's. Using this, claim 2 allows to determine
which sets of vertices form the vertices of some $k$-dimensional face. $\square$

The proof gives a poor algorithm 
and it was an interesting problem
to find better algorithms. This is an example where seeking an efficient 
algorithm was not motivated by questions from 
computer science but rather a natural aspect of our mathematical understanding.
Friedman (2009) found a remarkable LP-based polynomial time algorithm to
tell a simple polytope from its graph.
Another important open problem is to extend the theorem to dual graphs of arbitrary
triangulations of $(d-1)$-dimensional spheres. This is related to deep connections between polytopes and
spheres and various areas in commutative algebra and 
algebraic geometry, pioneered by Richard Stanley.

\end {document}